\documentclass[13pt]{amsart}

\usepackage{latexsym}
\usepackage{amsmath, amsfonts, amssymb}
\usepackage{graphicx,epsf}
\usepackage{enumerate}

\newtheorem{theorem}{Theorem}[section]

\newtheorem{corollary}{Corollary}[section]

\newtheorem{example}{Example}[section]
\newtheorem{algorithm}{Algorithm}[section]

\newtheorem{remark}{Remark}[section]

\numberwithin{equation}{section} \numberwithin{table}{section}
\numberwithin{figure}{section}

\newcommand{\Null}[1]{\mathcal{N}(#1)}
\renewcommand{\ker}[1]{\Null}

\everymath{\displaystyle}
\newcommand{\envert}[1]{\left\lvert#1\right\rvert}

\newcommand{\enVert}[1]{\left\lVert#1\right\rVert}

\newcommand{\eproof}{$\hfill\square$}

\begin{document}

\title[Adaptive ELM for Convection-Diffusion]{On Adaptive Eulerian--Lagrangian Method for \\ Linear Convection-Diffusion Problems}

\author[X. Hu]{Xiaozhe Hu}
\address{Department of Mathematics, The Pennsylvania State University, University Park, PA 16802}
\email{hu\_x@math.psu.edu}

\author[Y.-J Lee]{Youngju Lee}
\address{Department of the Mathematics, Rutgers, The State University of New Jersey, NJ 08854}

\author[J. Xu]{Jinchao Xu}
\address{Department of Mathematics, The Pennsylvania State University, University
Park, PA 16802} 
\thanks{}

\author[C.-S. Zhang]{Chensong Zhang}
\address{Institute of Computational Mathematics, Chinese Academy of Sciences, Beijing 100190, China.} 
\thanks{}

\date{}
\thanks{This work was supported in part by NSF Grant DMS-0915153, and DOE Grant DE-SC0006903.}

\begin{abstract}
In this paper, we consider the adaptive Eulerian--Lagrangian method (ELM) for linear convection-diffusion problems. Unlike the classical a posteriori error estimations, we estimate the temporal error along the characteristics and derive a new a posteriori error bound for ELM semi-discretization. With the help of this proposed error bound, we are able to show the optimal convergence rate of ELM for solutions with minimal regularity. Furthermore, by combining this error bound with a standard residual-type estimator for the spatial error, we obtain a posteriori error estimators for a fully discrete scheme. We present numerical tests to demonstrate the efficiency and robustness of our adaptive algorithm. 
\end{abstract}

\maketitle

\tableofcontents

\section{Introduction}\label{sec:intro}

Convection-diffusion equations oftentimes arise in mathematical models for fluid dynamics, environmental modeling, petroleum reservoir simulation, and other applications. Among them, the most challenging case for numerical simulation is the convection-dominated problems (when diffusion effect is very small compared with the convection effect)~\cite{Roos2008}. Dominated convection phenomena could appear in many real world problems; for example, convective heat transport with large P\'{e}clet numbers~\cite{Jakob1959}, simulation of oil extraction from underground reservoirs~\cite{Ewing1983}, reactive solute transport in porous media~\cite{Sorbents1996}, etc. The solutions of these problems usually have sharp moving fronts and complex structures; their nearly hyperbolic nature presents serious mathematical and numerical difficulties. Classical numerical methods, developed for diffusion-dominated problems,  suffer from spurious oscillations for convection-dominated problems. Many innovative ideas, like Upwinding, Method of Characteristics, and Local Discontinuous Galerkin methods, have been introduced to handle these numerical difficulties efficiently; see, for example, \cite{Abbott1966,Davis1976,Brooks1982a,Cockburn1998a} and references therein. 

For problems with nearly hyperbolic nature, it is nature to explore the idea of the so-called Method of Characteristics~\cite{Abbott1966}; and, this idea has been combined with different spatial discretizations like finite difference (FD), finite element (FE), and finite volume (FV) methods. Along this line of research, the Semi-Lagrangian method (or,  in the finite element context, the Eulerian--Lagrangian method) treats the convection and capacity terms together to carry out the temporal discretization in the Lagrangian coordinate. The Eulerian--Lagrangian method (ELM) gives rise to symmetric discrete linear systems, stabilizes the numerical approximation, and the corresponding diffusion problems are solved on a fixed mesh~\cite{Douglas1982,Pironneau1982}. This method and its variants have been successfully applied not only on the linear convection-diffusion problem~\cite{Douglas1982,Wang2000,Wang2007}, but also the incompressible Naiver-Stokes equations and viscoelastic flow problems; see, for example, \cite{Pironneau1982,Phillips1999,Phillips2000,Achdou2000,Xiu2001,Lee2006,Etienne2006}.

Adaptive mesh refinement (AMR) for partial differential equations (PDEs) has been the object of intense study for more than three decades. AMR techniques have been proved to be successful to deal with multiscale phenomena and to reduce the computation work without losing accuracy when solution is not smooth. In general, the adaptive algorithm for static problems generates graded meshes and iterations in the following form:
\begin{equation*}
\text{\textbf{SOLVE}} 
\rightarrow \text{\textbf{ESTIMATE}}
\rightarrow \text{\textbf{MARK}} 
\rightarrow \text{\textbf{REFINE/COARSEN}}
\end{equation*}
In the ESTIMATE procedure, we usually employ some computable local indicators to estimate the local error of the approximate solution we obtain from the SOLVE procedure. These indicators only depend on the datum and/or the approximate solution, and show in which part(s) of the domain the local error is relatively too big or too small in the current mesh. We then MARK these
subdomains and REFINE or COARSEN them accordingly.

Local error indicators determine whether the whole adaptive procedure is effective or not. Therefore, a reliable and efficient error indicator is the key to a successful AMR method and a posteriori error analysis is often used to obtain such an error indicator in practice~\cite{Ainsworth2000, Verfurth1996}. In the context of finite element methods, the theory of a posteriori error analysis and adaptive algorithms for linear elliptic problem is now rather mature. Convergence and optimality of adaptive methods for linear elliptic problems have been proved as the outcome of a sequence of work~\cite{Dorfler1996,Morin2002,Binev2004a,Stevenson2005,Cascon}; see the recent review by Nochetto, Siebert, and Veeser~\cite{Nochetto2009} and references therein.  On the other hand, for the nonlinear and time-dependent problems, the theory is still far from satisfactory. A posteriori error analysis for nonlinear evolution problems is even more challenging.

Adaptivity time-stepping is very important for time dependent problems because the practical problems sometimes have singularities or are multiscale in time. Uniform time step size cannot capture these phenomena. There are considerable amount of work in the literature devoted to the development of efficient adaptive algorithms for evolution problems. A posteriori error estimators for linear parabolic problems are studied in~\cite{Johnson1990,Eriksson1991,Eriksson1995,Picasso1998,Adjerid1999,Chen2004a,Lakkis2006} and are also derived for nonlinear problems; see~\cite{Moore1994,Eriksson1995a,Verfurth1998,Kruger2003} for example. There have been also some efforts for extending a posteriori error analysis for the time-dependent Stokes as well as the Navier-Stokes equations~\cite{Bernardi2004,Kharrat2005,Karakatsani2007}. In particular, a posteriori error estimates for convection-diffusion problems have been discussed in~\cite{Hebeker1999,Cawood2000,VERFURTH2005}.

It is nature to employ ARM techniques for convection-dominated problems because of the complex structures of the solutions and evolution of these structures in time. We also notice that spatial mesh adaptivity plays an important role in ELM to reduce numerical oscillations and smearing effect when inexact numerical integrations are employed~\cite{Jia2011}. Early adoption of adaptive characteristic methods has been seen since late 1980's~\cite{Demkowicz1986,Yeh1995,Chen2000}. A posteriori error estimates for characteristic-Galerkin method for convection-dominated problems have been proposed: Demokowicz and Oden~\cite{Demkowicz1986} considered the method of characteristics combined with Petrov-Galerkin finite element method for spatial discretization. Houston and S\"uli~\cite{Houston2001} give an a posteriori error estimator in $L^2(L^2)$-norm for the linear time-dependent convection-diffusion problem using the duality argument. Chen and Ji~\cite{Chen2004, Chen2006} give sharp a posteriori error estimations of ELM in $L^\infty(L^1)$-norm for linear and nonlinear convection-diffusion problems, respectively. A related $L^1(L^1)$ a posteriori error bound can be found in Chen, Nochetto, and Schmidt~\cite{Chen2000} for the continuous casting problem (convection-dominated and nonlinearly degenerated). 

In the previous error estimators mentioned above, the time residual, which measures the difference between the solutions of two successive time steps, is employed as a local indicator for temporal error. On the other hand, it is observed that ELM is essentially transforming a convection-diffusion problem to a parabolic problem along the characteristics. In this paper, we consider a posteriori error estimators for ELM, with focus on temporal error estimators. Motivated by Nochetto, Savar\'{e}, and Verdi~\cite{Nochetto2000}, in which adaptive time-stepping scheme for an abstract evolution equation $\partial_t u + \mathfrak{F}(u) = 0$ in Hilbert spaces is analyzed, we can obtain an a posteriori error bound for the temporal error along the characteristics. Combined with space adaptivity, adaptive method has been designed and implemented. Our numerical experiments in Section~\ref{sec:numer} suggest that the new a posteriori error estimator is effective. Moreover, the numerical results also indicate that the proposed error estimators are very efficient and can take advantage of the fact that ELM allows larger time stepsize.

The outline of this paper is as follows. In Section~\ref{sec:ELM}, we introduce a model problem and its Eulerian--Lagrangian discretization. In Section~\ref{sec:time}, we discuss temporal a posteriori and a priori error analysis for the linear convection-diffusion problem. In Section~\ref{sec:full}, we present a posteriori error estimation for the full-discretization. To simplify the discussion, we only consider a linear convection-diffusion model problem and restrict ourselves to the standard residual-type estimator for the spatial error, although the technique discussed in this paper can be potentially extended to other problems and spatial error estimators. In Section~\ref{sec:numer}, we consider the implementation of our adaptive algorithm and present some numerical experiments.

Throughout this paper, we will use the following notation. The symbol $L^2$ denotes the space of all square integrable functions and its norm is denoted by $\|\cdot\|$. Let $H^k$ be the standard Sobolev space of the scalar function whose weak derivatives up to order $k$ are square integrable, and, let $\|\cdot\|_k$ and $|\cdot|_k$ denote the standard Sobolev norm and its corresponding seminorm on $H^k$, respectively. Furthermore, $\|\cdot\|_{k,\omega}$ and $|\cdot |_{k,\omega}$ denote the norm $\|\cdot\|_k$ and the semi-norm $|\cdot |_k$ restricted to the domain $\omega\subset\Omega$, respectively. We also use the notation $H^k_0$ for the functions that belong to $H^k$ and their trace vanish on $\partial \Omega$. The dual space of $H^k_{0}$ is denoted by $H^{-k}$. We use the notation $X\lesssim(\gtrsim)Y$ to denote  the existence of a generic constant $C$, which depends only on $\Omega$, such that $X\leq (\geq) CY$.

\section{A model problem and discretization}\label{sec:ELM}

In this paper, we consider the following linear convection-diffusion model problem
\begin{equation} \label{eqn:model}
\frac{\partial u}{\partial t} + b(x,t) \cdot \nabla u - \varepsilon \Delta u = f  \qquad
\text{in} \ \Omega, \quad \text{for} \ t \in (0, T]
\end{equation}
with the initial condition 
\begin{equation}\label{eqn:ini_cond}
u(0) = u_0  \qquad \text{in} \ \Omega
\end{equation} 
and the boundary condition 
\begin{equation}\label{eqn:bd_cond}
u = 0 \qquad \text{on}  \ \partial \Omega, \quad \text{for} \ t \in (0,T]. 
\end{equation}
Here $\Omega \subset \mathbb{R}^d$ ($d=1,2,3$) is bounded polygonal domain.  We assume that $b$ is divergence free and vanishes on $\partial \Omega$ for $t \in (0, T]$. We assume that $u_0 \in L^2(\Omega)$ and $f \in L^2(0,T; H^{-1}(\Omega))$.

\subsection{The Eulerian--Lagrangian method}
In this section, we briefly recall the construction of the Eulerian--Lagrangian method. As usual, we introduce the characteristics (particle trajectory) $x(t; s, X)$, where $X$ is the \emph{Lagrangian coordinate} (original labeling) at time $s$ of the particle. It is also referred to the material coordinate. And $x$ is the \emph{Eulerian coordinate} at the current time and referred to as reference coordinate. Then for the given velocity field $b(x,t)$, the characteristics is defind by the following ordinary differential equation:
\begin{equation*}
\frac{\mathrm{d} x(t; s, X)}{\mathrm{d}t} = b(x(t;s,X),t), \quad x(s;s,X) = X
\end{equation*}
In the Eulerian coordinates, be chain rule, the material derivative is defined as following
\begin{equation}  \label{eqn:mate_deri}
\frac{\mathrm{D}u}{\mathrm{D}t} := \frac{\mathrm{d} u(x(t;s,X),t)}{\mathrm{d}t} = \frac{\partial u}{\partial t} + b \cdot \nabla u
\end{equation}
and then we can rewrite \eqref{eqn:model} as follow
\begin{equation}\label{eqn:SL_model}
\frac{\mathrm{D} u}{\mathrm{D} t} - \varepsilon \Delta u= f.
\end{equation}

In order to avoid the deformation of the mesh, ELM uses a fixed spatial mesh at each time-step and traces back along the characteristics. The characteristics at each time interval $[t_{n-1}, t_n]$ ($n=1,2,\ldots,N$) has the same original labeling $X$. Assume that $x^n(t) := x(t; t_n, X)$ and  
\begin{equation}\label{eqn:chara_ode}
\frac{\mathrm{d} x^n(t)}{\mathrm{d}t} = b(x^n(t),t), \quad x^n(t_n) = X, \quad t \in [t_{n-1}, t_n].
\end{equation}
We discretize the material derivative in \eqref{eqn:SL_model} using the backward Euler method as follows
\begin{equation}\label{eqn:TD_SL_model}
\frac{U^{n} - U^{n-1}(x^n(t_{n-1}))}{k_{n}} - \varepsilon \Delta U^{n} = f^{n}.
\end{equation}
Here, $k_n := t_{n} - t_{n-1}$, $f^n := f(t^n)$, and $U^n := U^n(x^n(t_n)) = U^n(X)$. It is easy to see that $x^n(t_{n-1}) = x^n(t_{n-1};t_n, X )$ is a function of $X$ in the Lagrangian coordinate. We usually refer to \eqref{eqn:TD_SL_model} as the temporal semi-discretization scheme. 

\subsection{Feet searching}
Our a posteriori error estimate base on the assumption that the characteristics are solved exactly, which preserves the determinant of the Jacobian of the flow. In the numerical analysis, it is difficult to do that and it was apparent that a computational realization of preserving the determinant of the Jacobian for the flow map to be one was crucial. Therefore, we discuss how to integrate the following
ordinary differential equation for the computation of the characteristic feet.  The numerical scheme we discuss here has second order accuracy and preserves the determinant of the Jacobian of the flow. 

\begin{eqnarray}\label{EQ:SOURCE}
\frac{d}{ds} y(x,t,s) &=& b(y(x,t,s),s),\\
y(x,t,t) &=& x.\nonumber
\end{eqnarray}
where
\begin{equation}\label{DIV:Free}
\nabla \cdot b = 0.
\end{equation}

The equation (\ref{EQ:SOURCE}) is often called the source-free
dynamical systems due to (\ref{DIV:Free}). For such a system, the
solution $y(x,t,s) : \mathbb{R}^{d} \mapsto \mathbb{R}^{d}$ is often called
the flow map or phase flow and it has the following property.

\begin{equation}\label{DIV:FreeII}
{\rm det} \left ( \frac{\partial y(x,t,s)}{\partial x} \right ) =
1, \quad \forall x \in \mathbb{R}^{d}, \,\, \forall s \in \mathbb{R}.
\end{equation}

We shall begin with introducing some popular scheme to
compute (\ref{EQ:SOURCE}) and showing that the scheme is indeed volume-preserving scheme for $d=2$ but not for $d=3$. We will then introduce some volume
preserving scheme to solve (\ref{EQ:SOURCE}) for $d=3$,
which is due to Feng and Shang~\cite{Feng1995}. 

In literatures, the following second order numerical scheme for
solving (\ref{EQ:SOURCE}) is popular and it seems to first
appear in \cite{CTemperton_AStaniforth_1987}.

First, we integrate the equation (\ref{EQ:SOURCE}) using the
mid-point rule to obtain :

\begin{equation}\label{APP:SOURCE}
\frac{1}{k}\left ( x - y(x,t,s) \right ) = \frac{1}{2} \, b
\left ( y \left (x,t,s+\frac{k}{2}\right ),s+\frac{k}{2} \right )
+ O(k^2),
\end{equation}
where $k = t-s$.

Second, the right hand side is approximated by a second order
accurate extrapolation. Namely,

\begin{equation}\label{EQ:EXTRAP}
b \left ( x,s+\frac{k}{2} \right ) = \frac{3}{2} b(x,s) -
\frac{1}{2} b (x,s-k) + O(k^2).
\end{equation}

The following approximation shall also be used :
\begin{equation}
y \left ( x,t,s+\frac{k}{2}\right ) = \frac{x + y(x,t,s)}{2} +
O(k^2).
\end{equation}

For notational conveniences, let us denote $y^{n}=y(x,t,s)$ and
$y^{n+\frac{1}{2}} = (x+y^n)/2.$ Hence we have the following
implicit approximations :
\begin{equation}\label{Mid:Rule}
\frac{1}{k} (x - y^n) := \frac{1}{2} \left ( \frac{3}{2}
b(y^{n+\frac{1}{2}},s) - \frac{1}{2} b(y^{n+\frac{1}{2}},s-k)
\right ).
\end{equation}

To see that the mid-point rule (\ref{Mid:Rule}) results in the
volume-preserving scheme, let us take the derivative with respect
to $x$ for both sides of  (\ref{Mid:Rule}). We then obtain the
following :
\begin{equation}\label{EQ:MID}
\frac{1}{k} \left (I - F (s) \right ) = A(y^{n+\frac{1}{2}})
\left ( I + F(s) \right ),
\end{equation}
where
\begin{equation}
A(y^{n+\frac{1}{2}}) = \frac{1}{2} \left ( \frac{3}{2} \nabla
b (y^{n+\frac{1}{2}},s) - \frac{1}{2} \nabla
b (y^{n+\frac{1}{2}},s-k) \right )
\end{equation}
and
\begin{equation}
F(s) = \frac{\partial y^n}{\partial x}.
\end{equation}

We solve (\ref{EQ:MID}) for $F(s)$ and obtain that
\begin{equation}
F(s) = \left ( I + \frac{k}{2} A(y^{n+\frac{1}{2}}) \right )^{-1}
\left ( I - \frac{k}{2} A(y^{n+\frac{1}{2}}) \right )
\end{equation}

Under the assumption that some appropriate finite element space
for the velocity field is used so that the divergence free
condition of the velocity field is imposed in the discrete sense,
we have
\begin{equation}\label{EQ:TR}
{\rm tr}A(y^{n+\frac{1}{2}}) = 0.
\end{equation}
From (\ref{EQ:TR}), we conclude that ${\rm det} F(s) = 1$
identically.

The main reason why such an algorithm is popular seems that it has
the volume preserving property. On the other hand, it is easy to
see that the algorithm may not result in the volume-preserving
scheme for $d=3$. Note that for $d=3$, under the assumption that
${\rm tr}A = 0$, for $H$ given as follows,
\begin{equation}
H = \left ( I + \frac{k}{2} A \right )^{-1} \left ( I -
\frac{k}{2} A \right )
\end{equation}
we have that 
$$
{\rm det} H = 1 \quad \Leftrightarrow  \quad {\rm det}A = 0.
$$ 

Our purpose here is that by reviewing the volume-preserving scheme
in three dimension developed by Feng and Shang in
\cite{Feng1995}, we wish to make sure such a volume
preserving scheme can be devised in three dimension and confirm
our numerical scheme can be implemented. As far as the author is
concerned, such a special algorithmic detail has not been
implemented in the context of the semi-Lagrangian scheme.

The basic idea of constructing the volume preserving scheme for $d
= 3$ is based upon the following observation. Following the idea
of H. Weyl, we have :

\begin{equation}\label{HW:Decomp}
b(y(x,t,s),s) = \left ( \begin{array}{cc} 0 \\ \frac{\partial
v_1}{\partial y_3}
\\ -\frac{\partial v_1}{\partial y_2} \end{array} \right ) +
\left ( \begin{array}{cc} \frac{\partial v_2}{\partial y_2}
\\ -\frac{\partial v_2}{\partial y_1} \\ 0 \end{array} \right ),
\end{equation}
where
\begin{eqnarray}
v_1 &=& -\int_{y_3}^{x_3} \left ( b_2(y_1,y_2,w,s) +
\frac{\partial b_2}{\partial y_1}(y_1,y_2,w,s) \right )\, dw \\
&+& \int_{y_2}^{x_2}
b_3(y_1,w,x_3,s)\,dw  \quad \mbox{ and } \\
\\
v_2 &=& -\int^{x_2}_{y_2} b_1(y_1,w,y_3,s)\,dw \nonumber
\end{eqnarray}

The actual expressions for $\frac{\partial v_1}{\partial y_2}$,
$\frac{\partial v_1}{\partial y_3}$, $\frac{\partial v_2}{\partial
y_2}$ and $\frac{\partial v_2}{\partial y_1}$ as follows :
\begin{eqnarray}
\frac{\partial v_1}{\partial y_2}                                                           &=& - b_3(y_1,y_2,y_3,s) \\
\frac{\partial v_1}{\partial y_3} &=& b_2(y_1,y_2,y_3,s)  - \int^{x_2}_{y_2} \frac{\partial b_1}{\partial y_1} (y_1,w,y_3,s)\,dw  \\
\frac{\partial v_2}{\partial y_2} &=& b_1(y_1,y_2,y_3,s) \\
\frac{\partial v_2}{\partial y_1} &=& -\int^{x_2}_{y_2}
\frac{\partial b_1}{\partial y_1} (y_1,w,y_3,s)\,dw.
\end{eqnarray}
From this, it is easy to see that (\ref{HW:Decomp}) holds and by
construction $u^1$ and $u^2$ given as follows are divergence
free :
\begin{eqnarray}
b^1 &=& \left ( \begin{array}{c} 0 \\ \frac{\partial
v_1}{\partial y_3}
\\ -\frac{\partial v_1}{\partial y_2} \end{array} \right )
 \quad \mbox{ and }
\quad  \\
b^2 &=& \left ( \begin{array}{c} \frac{\partial v_2}{\partial
y_2}
\\ -\frac{\partial v_2}{\partial y_1}
\\ 0 \end{array} \right )
\end{eqnarray}

Let us now denote $S_{i}^k$ by the volume preserving scheme for
\begin{eqnarray}\label{EQ:SFREE}
\frac{d}{ds}y(x,t,s) &=& b^i(y(x,t,s),s), \\
y(x,t,s) &=& x
\end{eqnarray}

 with $i=1,2$, then the following composition
is trivially volume preserving :
\begin{equation}\label{Com:Vol}
y^n = S_{2}^k \circ S_{1}^k
\end{equation}

Moreover, assuming $S_i^k$ is of second order accurate, it is easy
to see that the composition (\ref{Com:Vol}) is of second order.
The above idea on the composition is from Feng and Shang,
\cite{Feng1995}.

\subsection{Full discretization}
Discretizing \eqref{eqn:TD_SL_model} with suitable spatial discretizations, we obtain  fully discrete numerical schemes. In this paper, we focus on the finite element method. First, we define the weak forms of \eqref{eqn:SL_model} and \eqref{eqn:TD_SL_model} as usual. We define the bilinear form $a(\cdot, \cdot)$ as follow
\begin{equation}\label{def:bilinear_a}
a(v,w) := \varepsilon (\nabla v, \nabla w) \qquad \forall v,w \in V = H^1_0(\Omega).
\end{equation}
We denote the potential by
\begin{equation}\label{def:phi}
\phi(w) := \frac{1}{2}a(w,w) = \frac{\varepsilon}{2}\int_\Omega \envert{\nabla w}^2 \mathrm{d}x
\end{equation}
and its Frechet derivative by $\mathfrak{F}$, i.e.,
\begin{equation}\label{def:f}
(\mathfrak{F}(w), v) := a(w,v).
\end{equation}

It is easy to see that $\phi$ is convex and $\mathfrak{F}$ satisfies the following
inequality
\begin{equation}\label{ine:convex}
(\mathfrak{F}(w), v-w) \le \phi(v) - \phi(w) \qquad \forall w, v \in V.
\end{equation}
In fact, it is well-known that we have the following identity, 
\begin{equation} \label{eqn:identity-a}
a(w,w-v) = \phi(w) - \phi(v) + \frac{1}{2}|\!|\!|w-v|\!|\!|^{2} \qquad \forall w,v \in V,
\end{equation}
where $|\!|\!|v|\!|\!|^{2} := a(v,v)$. Furthermore, by taking the test function as $u-v$ in \eqref{eqn:SL_weak} and applying \eqref{eqn:identity-a}, we obtain that
\begin{equation}\label{eqn:VI}
(\frac{Du}{Dt}, u-v) + \phi(u) - \phi(v) + \frac{1}{2} |\!|\!| u-v |\!|\!|^{2} = 0.
\end{equation}

The weak forms of (\ref{eqn:SL_model}) and (\ref{eqn:TD_SL_model}) can be written as
\begin{equation}\label{eqn:SL_weak}
(\frac{\mathrm{D}u}{\mathrm{D}t}, v) + a(u,v) = (f,v) \qquad \forall v \in V,
\end{equation}
and 
\begin{equation}\label{eqn:TD_SL_weak}
(\frac{U^n-U^{n-1}(x^n(t_{n-1}))}{k_n}, v) + a(U^n, v) = (f^n,v) \qquad \forall v \in V,
\end{equation}
respectively. 

On a shape-regular triangulation $\mathcal{T}_{h}:=\{\tau_i\}$ of
$\Omega$, we introduce the piecewise continuous linear finite element space $V_h
\subset V$ such that
\begin{equation} \label{def:V_h}
V_h := \{v_h \in C(\overline{\Omega}): v_h | _{\tau_i} \in P^1(\tau_i),
\forall \tau_i \in \mathcal{T}_h \} \cap V.
\end{equation}
Let $\mathcal{T}_h^n$ and $V_h^n$ denote the mesh and the finite element space at time $t_n$ respectively. Then the fully discrete scheme can be written as: Suppose $U^{n-1}_h \in V_h^{n-1}$ is known,
find $U^n_h \in V_h^n$ such that
\begin{equation} \label{eqn:FD_SL_model}
(\frac{U^{n}_h - U^{n-1}_h ( x^n(t_{n-1}))}{k_{n}}, v_h) +
a (U^{n}_h, v_h) = (f_h^n,v) \qquad
\forall v_h \in V_h^n.
\end{equation}
where $f_h^n \in V_h^n$ is some suitable approximation of $f^n=f(t_n)$.

\section{Error analysis for temporal semi-discretization}\label{sec:time}

In this section, we focus on the a posteriori error estimation for temporal semi-discretization (\ref{eqn:TD_SL_model}). For the sake of simplicity, we assume that $f = 0$ in this section. Different from the standard time error indicators for ELM, which measure the solution difference along the time direction, our new error indicator measures the difference along the characteristics. We will establish a posteriori error estimation based on the new time error indicator and show its efficiency. We will also show an optimal priori error bound as a byproduct. Here optimality does not only mean the optimal convergence rate which can also be achieved by classic error analysis, but also mean the optimal regularity requirement which have not been proved by standard a priori error analysis. 

\subsection{A posteriori error analysis}
We first show that the solution of \eqref{eqn:SL_model} satisfies an energy identity along characteristics. In fact, by taking inner product with
 $\frac{\mathrm{D}u}{\mathrm{D}t}$ on both side of \eqref{eqn:SL_model} and applying the following identity 
\begin{equation*}
\frac{\mathrm{D}}{\mathrm{D}t} \phi(u) = (\mathfrak{F}(u), \frac{\mathrm{D}u}{\mathrm{D}t}),
\end{equation*}
we immediately obtain that
 \begin{equation}\label{eqn:energy}
\enVert{\frac{\mathrm{D} u}{\mathrm{D} t}}^{2} +
\frac{\mathrm{D}}{\mathrm{D}t} \phi(u) = 0.
\end{equation}
We can see that the convection diffusion equation preserve the total energy along  characteristics in continuous level. 

On the other hand, by choosing the test function $v
= (U^{n} - U^{n-1}(x^n(t_{n-1})))/k_{n} \in V$ in \eqref{eqn:TD_SL_weak} and employing \eqref{eqn:identity-a}, we have an discrete energy inequality: For any integer $1 \le n \le N$, 
\begin{equation}\label{ine:energy}
\left\|\frac{U^{n} - U^{n-1}(x^n(t_{n-1}))}{k_n}\right\|^2 +
\frac{\phi(U^{n}) - \phi(U^{n-1}(x^n(t_{n-1})))}{k_n} \leq 0.
\end{equation}
We note that the equality holds in \eqref{ine:energy} if there is no temporal discretization error. This motivates the following definition:
\begin{equation} \label{def:xi}
\xi_{n} := - \left\|\frac{U^{n} - U^{n-1}(x^n(t_{n-1}))}{k_n}\right\|^2 -
\frac{\phi(U^{n}) - \phi(U^{n-1}(x^n(t_{n-1})))}{k_n},
\end{equation}
and we can view the computable quantity $\xi_{n} \ge
0$ as a measure of the
deviation of numerical solution from satisfying the energy
conservation~\eqref{eqn:energy}. We can use $\xi_{n}$ as a
time error indicator for adaptive time stepping in ELM.

In order to give an a posteriori error bound for our new time error indicator, we construct the following linear interpolation 
\begin{equation}\label{eqn:linear_inter}
U(t) :=  \frac{t-t_{n-1}}{k_n} U^n(x^n(t_n)) + \frac{t_n - t}{k_n} U^{n-1} (x^n(t_{n-1})),
\end{equation}
where
\begin{equation*}
U^n(x^n(t_n)) = U^n(X) \quad \text{and} \quad U^{n-1}(x^n(t_{n-1})) = U^{n-1}(x^n(t_{n-1};t_n,X)).
\end{equation*}
Since the original labeling $X$ of the characteristics does not depend on time, we have
\begin{equation}\label{eqn:U_deri}
\frac{\mathrm{d} U(t)}{\mathrm{d} t} =  \frac{U^n - U^{n-1}(x^n(t_{n-1}))}{k_n} \qquad \forall t\in (t_{n-1},t_n].
\end{equation}
Substituting \eqref{eqn:U_deri} back to \eqref{eqn:TD_SL_weak} and then applying \eqref{eqn:identity-a}, we
have
\begin{equation*}
(\frac{\mathrm{d}U(t)}{\mathrm{d}t}, U^{n} - v) + \phi(U^{n}) - \phi(v)  + \frac{1}{2} |\!|\!| U^{n}-v |\!|\!|^{2} = 0 \qquad \forall v \in V.
\end{equation*}

By adding and subtracting $U(t)$, for any $v \in V$, we have 
\begin{equation}\label{ine:residual}
(\frac{\mathrm{d}U(t)}{\mathrm{d}t}, U(t) - v) + \phi(U(t)) - \phi(v) + \frac{1}{2} |\!|\!|U^{n}-v |\!|\!|^{2} =\mathfrak{R}(t),
\end{equation}
where $\mathfrak{R}(t)$ is the residual (which does not depend on the choice of test
function):
\begin{equation} 
\mathfrak{R}(t) := (\frac{\mathrm{d}U(t)}{\mathrm{d}t}, U(t) - U^{n}) + \phi(U(t)) - \phi(U^{n}).
\end{equation}
We notice that $\phi$ is a convexity functional, i.e.,
\begin{equation*}
\phi(U(t)) \le \frac{t_{n}- t}{k_n}\phi(U^{n-1}(x^n(t_{n-1}))) +
\frac{t-t_{n-1}}{k_n}\phi(U^{n}) \qquad \forall t \in[t_{n-1},t_n].
\end{equation*}
Hence, $\mathfrak{R}(t)$ can be bounded by
\begin{equation}\label{ine:residual_bound}
\mathfrak{R}(t) \le (t_{n} - t) \xi_{n}.
\end{equation}
By choosing $v = U(t)$ in \eqref{eqn:VI}, and $v = u(t)$ in
\eqref{ine:residual}, adding these two equations, and using \eqref{ine:residual_bound}, we end up with the following inequality
\begin{equation}\label{ine:error}
\frac{\mathrm{d}}{\mathrm{d}t} \enVert{u(x^n(t), t) - U(t)}^{2} + |\!|\!| u - U(t) |\!|\!|^{2} + |\!|\!| u - U^{n} |\!|\!|^{2}  \le
2 (t_{n} - t) \xi_{n} \quad t \in (t_{n-1}, t_{n}].
\end{equation}
Then the following upper bound of the new time error indicator holds:
\begin{theorem} \label{thm:TD_upper}
Let $u$ be the exact solution of \eqref{eqn:model} and
$\{U^{n}\}_{n=0}^{N}$ be the time semi-discrete solution in \eqref{eqn:TD_SL_model}. Assume that $\xi_n$ is defined as 
\eqref{def:xi}. For any integer $1 \le m \le N$ the folltowing upper bound holds:
\begin{equation} \label{ine:upper-bound}
\enVert{u(t_{m}) - U^{m}}^2  + \sum_{n=1}^{m} \int_{t_{n-1}}^{t_{n}} |\!|\!| u - U^{n} |\!|\!|^{2} \mathrm{d}t \le 
\sum_{n=1}^{m}k_n^2\xi_n.
\end{equation}
\end{theorem}

\proof  Integration \eqref{ine:error} in time, we can obtain
\begin{equation*}
\enVert{u(t_n) - U^n}^2  + \int_{t_{n-1}}^{t_{n}} |\!|\!| u - U^{n} |\!|\!|^{2} \mathrm{d}t  \le \enVert{u(x^n(t_{n-1}),t_{n-1}) - U^{n-1}(x^n(t_{n-1}))}^2 + k_n^2\xi_n.
\end{equation*}
Since $b$ is divergence free, which implies $\det \nabla x^n(t) = 1$, we have, after changing of variables, that
\begin{equation*}
\enVert{u(t_n) - U^n}^2 + \int_{t_{n-1}}^{t_{n}} |\!|\!| u - U^{n} |\!|\!|^{2} \mathrm{d}t \le \enVert{u(t_{n-1}) - U^{n-1}}^2 + k_n^2\xi_n.
\end{equation*}
Then \eqref{ine:upper-bound} follows directly from summing up
above inequality from $n=1$ to $m$. \eproof

\subsection{An optimal priori error estimation requiring minimal regularity}

Traditional a priori error analysis for ELM treats the temporal semi-discretization as a finite difference method and derive the error estimation based on the Taylor expansion (see \cite{Douglas1982,Pironneau1982}).  As a result, we obtain an optimal convergence rate but the regularity requirement is suboptimal. Taking advantage of the posteriori error estimation of new time error indicator $\xi$, we can derive an optimal order priori error estimation with minimal regularity requirement on the datum and the solution. 

From the definition of the new a posteriori error estimator \eqref{def:xi}, \eqref{eqn:TD_SL_weak} and \eqref{eqn:identity-a}, we have
\begin{equation*}
\xi_n = \frac{1}{2k_n} |\!|\!| U^n - U^{n-1}(x^n(t_{n-1}))|\!|\!|^{2}.
\end{equation*}
By the definition of $\mathfrak{F}(\cdot)$ and \eqref{eqn:TD_SL_model}, we have
\begin{equation*}
\xi_n = \frac{1}{2}(\mathfrak{F}(U^{n-1}(x^n(t_{n-1}))) - \mathfrak{F}(U^n),
\mathfrak{F}(U^n)),
\end{equation*}
Substitute back into \eqref{ine:upper-bound}, we have
\begin{equation*}
\max_{1 \le n \le N} \enVert{u(t^{n}) - U^{n}}^2 + \sum_{n} \int_{t_{n-1}}^{t_{n}} |\!|\!| u - U^{n} |\!|\!|^{2} \mathrm{d}t  \le 
\frac{1}{2}\sum_{n=1}^{N}k_n^2 (\mathfrak{F}(U^{n-1}(x^n(t_{n-1})) -
\mathfrak{F}(U^n), \mathfrak{F}(U^n)).
\end{equation*}
Set $k_{\max} := \max_{1 \le n \le N} k_n$ and use the elementary
inequality $2(v-w,v) \ge \envert{v}^2 - \envert{w}^2$, we obtain
\begin{multline} \label{ine:priori_error}
\max_{1 \le n \le N} \enVert{u(t^{n}) - U^{n}}^2 + \sum_{n} \int_{t_{n-1}}^{t_{n}} |\!|\!| u - U^{n} |\!|\!|^{2} \mathrm{d}t  \\
\le
\frac{1}{4}k_{\max}^2\sum_{n=1}^{N} (\enVert{\mathfrak{F}(U^{n-1}(x^n(t_{n-1})))}^2 -
\enVert{\mathfrak{F}(U^n)}^2).
\end{multline}

Now using the assumption that $b$ is divergence free and summing up \eqref{ine:priori_error} from $n=1$ to $m$, we obtain the following optimal a priori error estimation with minimal regularity requirement:
\begin{corollary}
Let $u$ be the solution of \eqref{eqn:model} with initial value $u_0$,  and $\{U^{n}\}_{n=1}^{N}$ is the time semi-discrete numerical solution of the temporal semi-discretization \eqref{eqn:TD_SL_model}. For any integer $1\le m \le N$, we have the following priori error estimate
\begin{equation}\label{ine:priori}
\|u(t^{m}) - U^{m}\|^2 + \sum_{n=1}^{m} \int_{t_{n-1}}^{t_{n}} |\!|\!| u - U^{n} |\!|\!|^{2} \mathrm{d}t  \le
\frac{1}{4}k_{\max}^2 \|\mathfrak{F}(u_0)\|^2.
\end{equation}
\end{corollary}

\begin{remark}[Regularity Requirement]
Only on initial guess.
\end{remark}

\section{A posteriori error estimation for full discretization}\label{sec:full}

In this section, a posteriori error estimation for the fully-discrete scheme \eqref{eqn:FD_SL_model} is considered. The characteristic feet $x^n(t_{n-1})$ is defined by \eqref{eqn:chara_ode} and is computed exactly. Besides the new time error indicator, residual type error estimator is used as a spatial error indicator. 
Let $R^n$ be the \emph{interior residual}, i.e.,
\begin{equation} \label{def:inter_residual}
R^n : = f^n_h - \frac{U^n_h - U^{n-1}_h(x^n(t_{n-1}))}{k_n} +
\varepsilon \Delta U^n_h, \quad \text{on any} \ \tau \in \mathcal{T}_h^n
\end{equation}
and $J_e^n$ be the \emph{jump residual} also defined, i.e.,
\begin{equation} \label{def:jump_residual}
J^n_e := ( \nabla U^n_h|_{\tau_1} - \nabla
U^{n}_h|_{\tau_2}) \cdot \nu_e, 
\end{equation}
where $e= \partial \tau_1 \cap
\partial \tau_2$ is the edge/face shared by $\tau_1 \in \mathcal{T}_h^n$ and $\tau_2 \in \mathcal{T}_h^n$, $\nu_e$ is the unit normal vector of $e$ from $\tau_1$ to $\tau_2$. 

As the analysis of time semi-discretization, we define the linear interpolation $U_h(t)$ of $\{U_h^n\}_{n=0}^N$ as following
\begin{equation}\label{def:FD_linear_inter}
U_h(t) = \frac{t-t_{n-1}}{k_n}U^n_h + \frac{t_n -t}{k_n}
U^{n-1}_h (x^n(t_{n-1})),
\end{equation}
Similar with \eqref{eqn:U_deri}, we have
\begin{equation}
\frac{\mathrm{d} U_h}{\mathrm{d}t} = \frac{U_h^n - U_h^{n-1} (x^n(t_{n-1}))}{k_n}.
\end{equation}
Now, the following a posteriori error estimation holds
\begin{theorem}\label{thm:FD_upper}
Let $u$ and $\{U_h^n\}_{n=0}^{N}$ be the solutions of
\eqref{eqn:model} and \eqref{eqn:FD_SL_model} respectively. For any integer $1\le m \le N$, there exists a constant $C >0$ which does not depend on
the mesh size and the time step size, such that the following a
posteriori error bound
\begin{equation}
\label{ine:FD_upper_bound}
\begin{split}
& \frac{1}{2}||u(t_m)-U_h^m||^2   + \frac{1}{2}\sum_{n=1}^{m} \int_{t_{n-1}}^{t_{n}} |\!|\!| u - U_{h}^{n} |\!|\!|^{2} \mathrm{d}t  \\
 \le & ||u(0) - U_h^0||^2+\sum_{n=1}^{m} k_n^2 \xi_n 
+ C\sum_{n=1}^{m} k_n \eta_n  
 + \sum_{n=1}^{m} \int_{t_{n-1}}^{t_n} ||f(x^n(t), t)- f^n_h||_{-a}^2 \mathrm{d}t
\end{split}
\end{equation}
holds, where the temporal error indicator $\xi_n$ and the spatial error indicator $\eta_n$ are defined as
\begin{align} \label{def:time_indicator}
 \xi _n := & \left( f_h^n -
\frac{U^n_h - U_h^{n-1}(x^n(t_{n-1}))}{k_n}, \frac{U_h^n -
U_h^{n-1}(x^n(t_{n-1}))}{k_n} \right)  \nonumber \\  &- \frac{\phi(U^n_h)
- \phi(U_h^{n-1}(x^n(t_{n-1})))}{k_n}
\end{align}
and 
\begin{equation}\label{def:space_indicator}
\eta_n := \sum_{\tau} \eta_{\tau}^n, \quad \eta^n_{\tau} :=
 \frac{1}{\varepsilon} h_{\tau}^2 \enVert{R^n}_{L^2(\tau)}^2 +
\varepsilon \sum_{e \in \tau} h_e \enVert{J^n_e}^2_{L^2(e)}.
\end{equation}
\end{theorem}

\proof
For any $ v \in V$ and $ v_h \in V_h^n$, we have
\begin{equation}
(\frac{\mathrm{d}U_h}{\mathrm{d}t} - f^n_h, v) + a(U^n_h, v) = (\frac{\mathrm{d} U_h}{\mathrm{d}t} -
f^n_h, v-v_h) + a(U^n_h, v-v_h).
\end{equation}
Let $w = U_h^n$ in \eqref{eqn:identity-a} and $v_h = U_h^n$ in \eqref{eqn:FD_SL_model}, we have
\begin{eqnarray*}
&(\frac{\mathrm{d} U_h}{\mathrm{d}t} - f^n_h, U^n_h -v) + \phi(U^n_h) -\phi(v) +
\frac{1}{2}|\!|\!|v-U_h^n|\!|\!|^2  \\ &= (f^n_h - \frac{ \mathrm{d} U_h}{\mathrm{d}t}, v-v_h)
- a(U^n_h, v-v_h).
\end{eqnarray*}
Inserting $U_h$, we have
\begin{eqnarray*}
&&(\frac{\mathrm{d} U_h}{\mathrm{d}t} - f_h, U_h - v) + \phi(U_h) - \phi(v) +
\frac{1}{2}|\!|\!|v-U_h^n|\!|\!|^2 \\ &=& (\frac{\mathrm{d} U_h}{\mathrm{d}t} - f_h, U_h -
U^n_h) + \phi(U_h) - \phi(U^n_h)  + (f^n_h - \frac{\mathrm{d} U_h}{\mathrm{d}t},
v-v_h) - a(U^n_h, v-v_h).
\end{eqnarray*}
Let $w=v$ and $v = U_h$ in \eqref{eqn:identity-a}, we have
\begin{eqnarray*}
&&(\frac{\mathrm{d} U_h}{\mathrm{d}t}, U_h - v) + a(v, U_h - v) + \frac{1}{2}|\!|\!|v -
U_h|\!|\!|^2 + \frac{1}{2}|\!|\!|v-U_h^n|\!|\!|^2 \cr &=& (\frac{\mathrm{d} U_h}{\mathrm{d}t} -
f_h, U_h - U^n_h) + \phi(U_h) - \phi(U^n_h) \\ && + (f^n_h -
\frac{\mathrm{d} U_h}{\mathrm{d}t}, v-v_h) - a(U^n_h, v-v_h) + (f_h^n, U_h -v).
\end{eqnarray*}
Then let $v = u $ and $v_h = U_h^n +\Pi_h^n(u - U_h^n) $ where $\Pi_h^n: V \rightarrow V_h^n$ is the so called Clement interpolation operator. Moreover, let $v = u-U_h $ in \eqref{eqn:SL_weak} and add to above inequality, we can obtain that
\begin{align*}
&\frac{1}{2}\frac{\mathrm{d}}{\mathrm{d}t}\enVert{u-U_h}^2 +
\frac{1}{2}|\!|\!|u - U_h|\!|\!|^2 + \frac{1}{2}|\!|\!|u-U_h^n|\!|\!|^2 \\ =&
(\frac{\mathrm{d} U_h}{\mathrm{d}t} - f^n_h, U_h - U^n_h) + \phi(U_h) - \phi(U^n_h)
\\  & +(f^n_h - \frac{\mathrm{d} U_h}{\mathrm{d}t}, u-U_h^n - \Pi_h^n(u-U_h^n)) - a(U^n_h, u-U_h^n - \Pi_h^n(u-U_h^n))  \\ &+ (f - f^n_h, u - U_h).
\end{align*}
Then by the definitions \eqref{def:inter_residual}, \eqref{def:jump_residual} and \eqref{def:time_indicator}, we have
\begin{align*}%
&\frac{1}{2}\frac{\mathrm{d}}{\mathrm{d}t} \enVert{u-U_h}^2 +
\frac{1}{2}|\!|\!|u-U_h|\!|\!|^2 + \frac{1}{2}|\!|\!|u-U_h^n|\!|\!|^2 \\ \le&
(t_n - t) \xi_n  \\& + (R^n, (u - U_h^n) - \Pi_h^n(u -
U_h^n)) + \sum_e \int_e J^n_e[(u - U_h^n) - \Pi_h^n(u -
U_h^n)]\mathrm{d}s \\ & +(f - f_h^n, u- U_h) \\
\le &(t - t_{n}) \xi_{n} + C \eta_{n}^{1/2} |\!|\!| u - U_{h}^{n} |\!|\!| +  |\!| f - f_h^n |\!|_{-a} |\!|\!| u- U_h |\!|\!|\\
\le &(t - t_{n}) \xi_{n} + C \eta_{n} + \frac{1}{4} |\!|\!| u - U_{h}^{n} |\!|\!|^{2} + \frac{1}{2}|\!|f - f^n_h|\!|^2_{-a} + \frac{1}{2} |\!|\!|u-U_h|\!|\!|^2.
\end{align*}
Here in last two inequalities, we use the interpolation property of the Clement interpolation,  the Cauchy Schwartz inequality, and the Young's inequality.  The norm $|\!| \cdot |\!|_{-a}$ is defined by $|\!| v |\!|_{-a} := \sup_{w} (v,w)/ |\!|\!| w |\!|\!|$.  Now, we have
\begin{equation}\label{ine:error-represent}
\frac{1}{2}\frac{\mathrm{d}}{\mathrm{d}t} \enVert{u(x^n(t),t)-U_h(t)}^2 +
\frac{1}{4}|\!|\!|u-U_h^{n}|\!|\!|^2 \le  (t_n-t)\xi_n 
+C \eta_n + \frac{1}{2}|\!|f - f^n_h|\!|^2_{-a}.
\end{equation}
For any $t^{*} \in (t_{m-1}, t_{m}]$, by integrating \eqref{ine:error-represent} in time from $t_{m-1}$ to $t_{m}$, we have
\begin{align*}
&\quad \frac{1}{2}||u(x(t^{*}),t^{*}) - U_h(t^{*})||^2  + \frac{1}{4}\int_{t_{m-1}}^{\min(t_{m}, t^{*})}|\!|\!|u-U_h^{n}|\!|\!|^2 \\ 
&\le \frac{1}{2}||u(x^m(t_{m-1}), t_{m-1}) - U_h^{m-1}(x^m(t_{m-1}))||^2  +\frac{1}{2} k_m^2 \xi_m   + C k_m \eta_{m} \\& + \frac{1}{2}\int_{t_{m-1}}^{\min(t_{m}, t^{*})} |\!| f (x^m(t),t) - f_h^m |\!|^2_{-a} \mathrm{d}t.
\end{align*}
For $n < m$, by integrating \eqref{ine:error-represent} in time from $t_{n-1}$ to $t_{n}$, we have
\begin{eqnarray*}
&\frac{1}{2}||u(t_n) - U_h^n||^2 + \frac{1}{4}\int_{t_{n-1}}^{t_n}|\!|\!|u-U_h^{n}|\!|\!|^2 \\ 
&\le \frac{1}{2}||u(x^n(t_{n-1}), t_{n-1}) - U_h^{n-1}(x^n(t_{n-1}))||^2  +\frac{1}{2} k_n^2 \xi_n  + C k_n
\eta_n + \frac{1}{2}\int_{t_{n-1}}^{t_n} |\!| f (x^n(t),t) - f_h^n |\!|_{-a}^2 \mathrm{d}t.
\end{eqnarray*}
Using the assumption that $b$ is divergence free, sum the above two inequalities from $n=1$ to $m$,  we have
\begin{align*}
&\quad \frac{1}{2}||u(x(t^{*}),t^{*}) - U_h(t^{*})||^2 + \frac{1}{4} \sum_{n=1}^{m} \int_{t_{n-1}}^{\min(t_{n}, t^{*})}|\!|\!|u-U_h^{n}|\!|\!|^2 \\ 
&\le \frac{1}{2}||u(0) - U_h^{0}||^2 +\frac{1}{2} \sum_{n=1}^{m} k_n^2 \xi_n  + C \sum_{n=1}^{m} k_n
\eta_{n}   + \frac{1}{2} \sum_{n=1}^{m} \int_{t_{n-1}}^{t_{n}} |\!| f (x^n(t),t) - f_h^n |\!|_{-a}^2 \mathrm{d}t.
\end{align*}
Note that $t^{*}$ is arbitrary, the above inequality implies the upper bound \eqref{ine:FD_upper_bound}.\eproof

\section{Numerical experiments}\label{sec:numer}

In this section, we first introduce the adaptive algorithm we use in numerical experiments and then present some numerical results.

\subsection{Adaptive Algorithm}
Since the main purpose of this paper is analyze the proposed time error estimator, we start with the algorithm for time step size control. Let $\mathtt{TOL_{time}}$ be the total tolerance allowed for the error introduced by the time discretization, according to Theorem \ref{thm:FD_upper}, that is
\begin{equation}\label{ine:time_control}
\sum_{n=1}^{N}k_n^2 \xi_n + 2(\sum_n \int_{t_{n-1}}^{t_n} ||f(x^n(t), t)- f^n_h||_{-a} \mathrm{d}t)^2 \le \mathtt{TOL_{time}}
\end{equation}
Therefore, a natural way to satisfy \eqref{ine:time_control} is to adjust the time step size such that
\begin{equation} \label{ine:time_ind_control}
k_n \xi_n \le \frac{\mathtt{TOL_{time}}}{2T}  \quad  \text{and} \quad  \frac{1}{k_n} \int_{t_{n-1}}^{t_n} \enVert{f(x^n(t),t) - f_h^n}_{-a} \mathrm{d}t \le \frac{1}{2T} \sqrt{\mathtt{TOL_{time}}}
\end{equation}
Based on this, now we propose the following algorithm to adjust the time step size and control the error under given tolerance $\mathtt{TOL_{time}}$.

\begin{algorithm}[H]\label{alm:time_adaptive}
Given tolerance $\mathtt{TOL_{time}}$, $\delta_1 \in (0,1)$, $\delta_2 > 1$, $\theta \in (0, 1)$, and the initial time stepsize $k_0$.
\begin{enumerate}
\item  Set $n:=1$
\item  Set $k_n := k_{n-1}$
\item  Solve the time semi-discretization problem using time step size $k_{n}$
\item  Compute the time error indicator \eqref{def:time_indicator}
\item  Check the error: If \eqref{ine:time_ind_control} is not satisfied, $k_n:= \delta_1 k_n$ and goto (2) 
\item  If $k_n \xi_n \le \theta \frac{\mathtt{TOL_{time}}}{2T} $ and $ \frac{1}{k_n} \int_{t_{n-1}}^{t_n} \enVert{f(x^n(t),t) - f_h^n}_{-a} \mathrm{d}t \le \frac{1}{2T} \sqrt{\theta \mathtt{TOL_{time}}}$, $k_n:= \delta_2 k_n$
\item  Let $n:=n+1$ and goto (2)
\end{enumerate}
\end{algorithm}

The above algorithm is for adjusting the time step size and controlling the error introduce by time discretiziation. When we consider the fully discretization, we need to take mesh adaption into account.
According to Theorem \ref{thm:FD_upper}, the space error indicator \eqref{def:space_indicator} gives us the elements where the local error is relatively large and should be refined. Let $\mathtt{TOL_{space}}$ be the tolerance allowed for the error introduced by the spatial discretization, then the usual stopping criterion for mesh refinement is:
\begin{equation} \label{ine:space_control}
\sum_n k_n \eta_n \le \mathtt{TOL_{space}}
\end{equation}
As usual, in order to achieve equal distribution of error, \eqref{ine:space_control} can be replaced by
\begin{equation}\label{ine:space_ind_control}
\eta_n \le \frac{\mathtt{TOL_{space}}}{T}. 
\end{equation}

The main difficulty for time-dependent problems is the mesh coarsening. We need a coarsening error indicator to guide the mesh coarsening procedure. Here, we employ the same coarsening error indicator as the one proposed in \cite{Chen2004a}. Let  $\mathcal{T}_H^n$ be a coarsening mesh from mesh $\mathcal{T}_h^n$, and $V_H^n$ and $V_h^n$ be the corresponding finite element spaces. It is nature that $V_H^n \subset V_h^n$. Furthermore, let $U_H^n \in V_H^n$ and $U_h^n$ be the numerical solution at time $t_n$ satisfy that
\begin{align*}
&(\frac{U_H^n - U_h^{n-1}(x^n(t_{n-1}))}{k_n}, v) + a(U_H^n, v) = (f^n, v), \quad \forall v \in V_H^n \\
&(\frac{U_h^n - U_h^{n-1}(x^n(t_{n-1}))}{k_n}, v) + a(U_h^n, v) = (f^n, v), \quad \forall v \in V_h^n
\end{align*}
Subtracting these two equalities and taking $v = U_H^n - I_H^n U_h^n \in V_H^n$, we have
\begin{equation*}
(U_{H}^{n} - U_{h}^{n}, U_{H}^{n} - I_{H}^{n}U_{h}^{n}) + k_{n} a(U_{H}^{n} - U_{h}^{n}, U_{H}^{n} - I_{H}^{n}U_{h}^{n}) = 0.
\end{equation*}
Using the elementary equality $2ab = a^{2} + b^{2} - (b-a)^{2}$ for the two terms on the left hand side of the above equality respectively, we can obtain the following Galerkin orthogonal relation
\begin{equation*}
\begin{split}
\enVert{U_H^n - U_h^n}^2 + k_n |\!|\!| U_H^n - U_h^n |\!|\!|^2  = & \enVert{U_h^n - I_H^n U_h^n}^2 + k_n |\!|\!|U_h^n - I_H^n U_h^n|\!|\!|^2 \\ -& \enVert{U_H^n - I_H^n U_h^n}^2 - k_n |\!|\!|U_H^n - I_H^n U_h^n|\!|\!|^2.
\end{split}
\end{equation*}
Then we have 
\begin{equation*}
\enVert{U_H^n - U_h^n}^2  \le \enVert{U_h^n - I_H^n U_h^n}^2 + k_n |\!|\!|U_h^n - I_H^n U_h^n|\!|\!|^2. 
\end{equation*}
Together with the triangular inequality
\begin{equation*}
\enVert{u(t_n) - U_H^n} \le \enVert{u(t_n) - U_h^n} + \enVert{U_H^n - U_h^n}
\end{equation*}
and Theorem \ref{thm:FD_upper}, this suggests us the following coarsening error indicator,
\begin{equation}\label{def:coarse_indicator}
\zeta_n :=  \sum_{\tau}\zeta^n_{\tau}, \quad \zeta_{\tau}^n := \enVert{U_h^n - I_H^n U_h^n}^2_{L^2(\tau)} + k_n |\!|\!|U_h^n - I_H^n U_h^n|\!|\!|^2_{L^2(\tau)},
\end{equation}
and the stopping criterion for coarsening
\begin{equation}\label{ine:coarsen_ind_control}
\zeta_n \le \frac{\mathtt{TOL_{coarsen}}}{T}.
\end{equation}
This indicator does not depend on the solution $U_H^n$, and allows us to coarsen only once after mesh refinement, which avoid the undesirable case that the elements which were marked for coarsening must be refined again to reduce the total error.

There are several different ways to couple time step size control and mesh adaption. A simple way is to include space adaptivity in every iteration of Algorithm \ref{alm:time_adaptive}. But this may result in solving a large number of discrete systems and adaptivity may suffer due to the computational overhead. Here we apply a mild way to couple the time and space adaptivity which is a modification of the adaptive procedure in \cite{Schmidt2005}. At each time step, we first adjust the time step size with the old mesh, and then adapt the mesh, and followed by checking the time error again.  
This leads to the following adaptive algorithm for on single time step.
\begin{algorithm}[H]
Given the tolerance $\mathtt{TOL_{time}}$, $\mathtt{TOL_{space}},$ and $\mathtt{TOL_{coarsen}}$, parameters $\delta_1 \in (0,1)$, $\delta_2 >1$ and $\theta \in (0,1)$. Let $U_h^{n-1}$ and the mesh $\mathcal{T}_h^{n-1}$ have been obtained at the previous time step $t_{n-1}$.
\begin{enumerate}

\item 
\begin{itemize}
\item Set $\mathcal{T}_h^{n} := \mathcal{T}_h^{n-1}$, $k_n:=k_{n-1}$ and $t_{n} = t_{n-1} + k_n$
\item Solve the fully discrete problem \eqref{eqn:FD_SL_model} on $\mathcal{T}_h^n$ using time step size $k_n$
\item Compute the error indicators
\end{itemize}

\item While \eqref{ine:time_ind_control} is not satisfied, 
\begin{itemize}
\item Set $k_n:= \delta_1 k_n$ and $t_{n} = t_{n-1} + k_n$
\item Solve the fully discretize problem \eqref{eqn:FD_SL_model} on $\mathcal{T}_h^n$ using time step size $k_n$
\item Compute the error indicators
\end{itemize}
End While 

\item While \eqref{ine:space_ind_control} is not satisfied,
\begin{itemize}
\item Refine the mesh and produce a refined $\mathcal{T}_h^n$
\item Solve the fully discretize problem \eqref{eqn:FD_SL_model} on $\mathcal{T}_h^n$ using time step size $k_n$
\item Compute the error indicators
\item While \eqref{ine:time_ind_control} is not satisfied, 
\begin{itemize}
\item Set $k_n:= \delta_1 k_n$ and $t_{n} = t_{n-1} + k_n$
\item Solve the fully discretize problem \eqref{eqn:FD_SL_model} on $\mathcal{T}_h^n$ using time step size $k_n$
\item Compute the error indicators
\end{itemize}
End While
\end{itemize}
End While 

\item 
\begin{itemize}
\item Coarsen the mesh to produce a new mesh $\mathcal{T}_h^n$ according to \eqref{ine:coarsen_ind_control}
\item Solve the fully discretize problem \eqref{eqn:FD_SL_model} on $\mathcal{T}_h^n$ using time step size $k_n$
\end{itemize}

\item  If  $k_n \xi_n \le \theta \frac{\mathtt{TOL_{time}}}{2T} $ and $ \frac{1}{k_n} \int_{t_{n-1}}^{t_n} \enVert{f(x^n(t),t) - f_h^n}_{-a} \mathrm{d}t \le \frac{1}{2T} \sqrt{\theta \mathtt{TOL_{time}}}$, then	
	\begin{itemize}
	\item $k_n:= \delta_2 k_n$
	\end{itemize}
	End If
\end{enumerate}
\end{algorithm}

\subsection{Numerical results}
Now we will present some numerical results to demonstrate the efficiency and robustness of our new adaptive method. In our numerical experiments, since we use backward Euler to discretize the material derivative, we choose $\delta_1 = 0.5$, $\delta_2 = 2$ and $\theta =0.5$. 

\subsubsection{1D Examples}
\begin{example} \label{exp:peak-1D}
Given the initial condition
\begin{equation*}
u(x, 0) = e^{-(x-x_{0})^{2}/2\lambda^{2}},
\end{equation*}
and $f=0$, the exact solution for the model problem~\eqref{eqn:model} (if $b$ is constant) is given by 
\begin{equation*}
u(x, t) = \frac{\lambda}{\sqrt{\lambda^{2} + 2\varepsilon
t}}e^{-(x-x_{0}-bt)^{2}/2(\lambda^{2} + 2\varepsilon t)},
\end{equation*}
where $\lambda$ is a parameter measures the width of the support of the solution. The computational domain is $[-1,2]$, and we give the Dirichlet boundary condition $u(-1) = u(2) = 0$.
\end{example}

\begin{example}\label{exp:shock-1D}
Given the initial condition
\begin{equation*}
u(x,0) = \left\{\begin{array}{ll}
1 & \quad x \leq 0\\
0 & \quad x > 0,
\end{array}
\right.
\end{equation*}
and $f=0$, the exact solution for the model problem~\eqref{eqn:model} (if $b$ is constant) is given by
\begin{equation*}
u(x,t) = \frac{1}{2}\left\{  \mathrm{erfc}\left(\frac{x-bt}{2\sqrt{\varepsilon t}}\right) + \exp\left(\frac{bx}{\varepsilon}\right)\mathrm{erfc}\left(\frac{x+bt}{2\sqrt{\varepsilon t}}\right) \right\},
\end{equation*}
where $\mathrm{erfc}(x) = \frac{2}{\sqrt{\pi}}\int^{\infty}_{x}e^{-s^{2}}\mathrm{d}s$ is the so-called complementary error function. The computational domain is $[0,2]$, and we give the Dirichlet boundary condition $u(0) = 1, u(2) = 0$.
\end{example}

\begin{figure}[ht]
\caption{Example \ref{exp:peak-1D}. Left: time step size comparison; Right: numerical solution comparison} \label{fig:peak-1D}
\centering
\includegraphics[width=0.49\textwidth, height=2.5in]{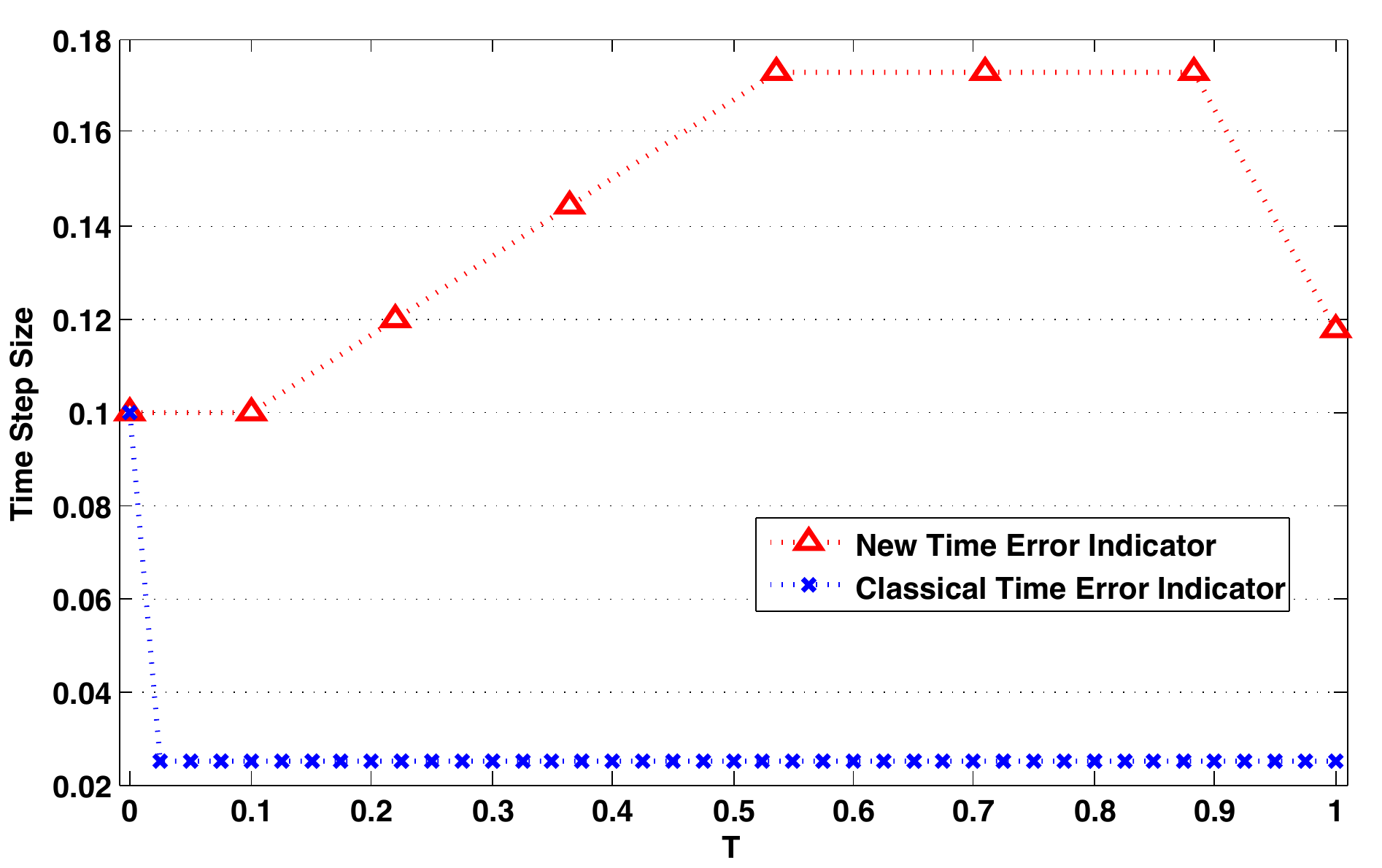}
\hfill
\includegraphics[width=0.49\textwidth, height=2.5in]{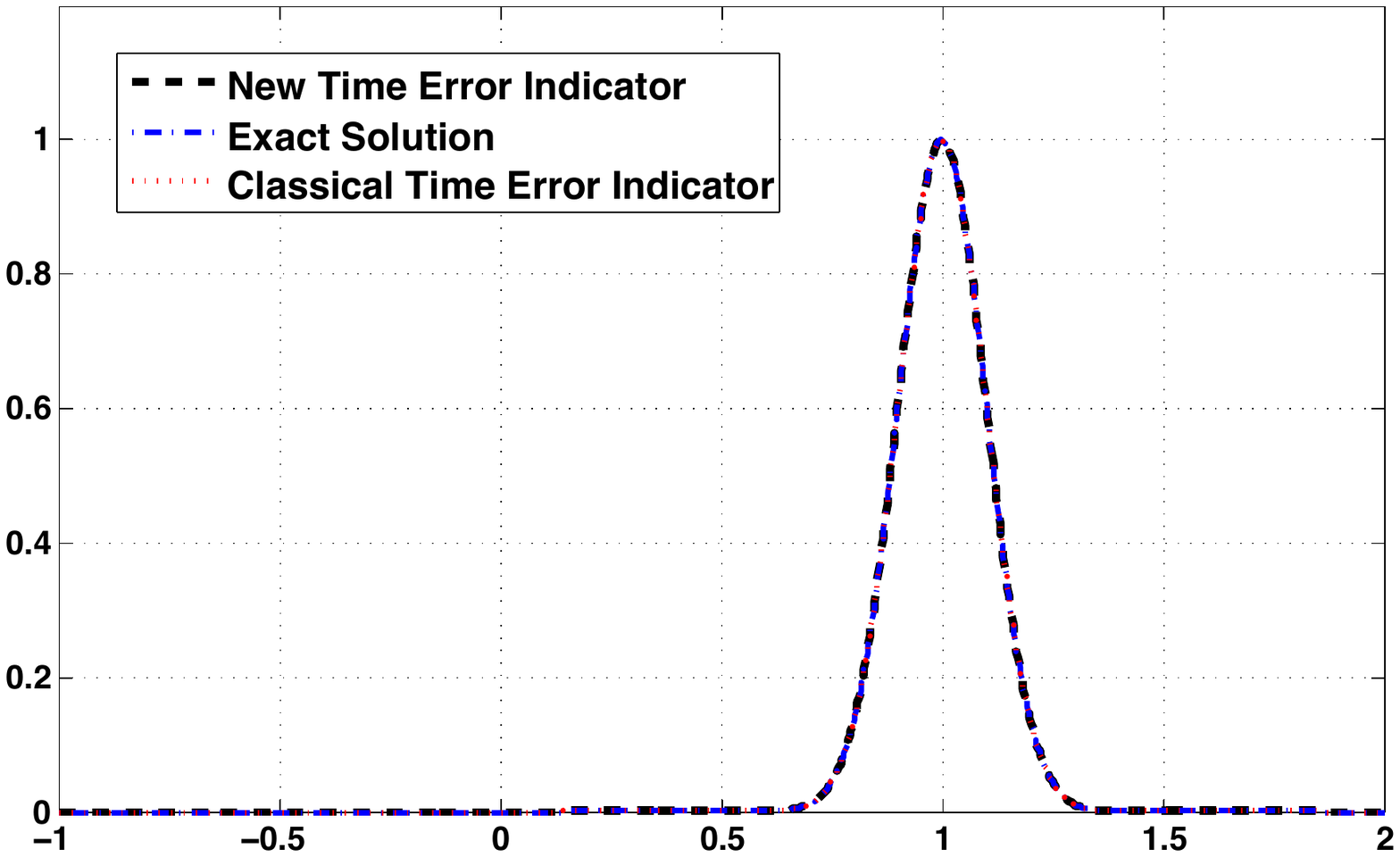}
\end{figure}
 
\begin{figure}[ht]
\caption{Example \ref{exp:shock-1D} ($\varepsilon = 10^{-6}$). Left: time step size comparison; Right: numerical solution comparison} \label{fig:shock-1D}
\centering
\includegraphics[width=0.49\textwidth, height=2.5in]{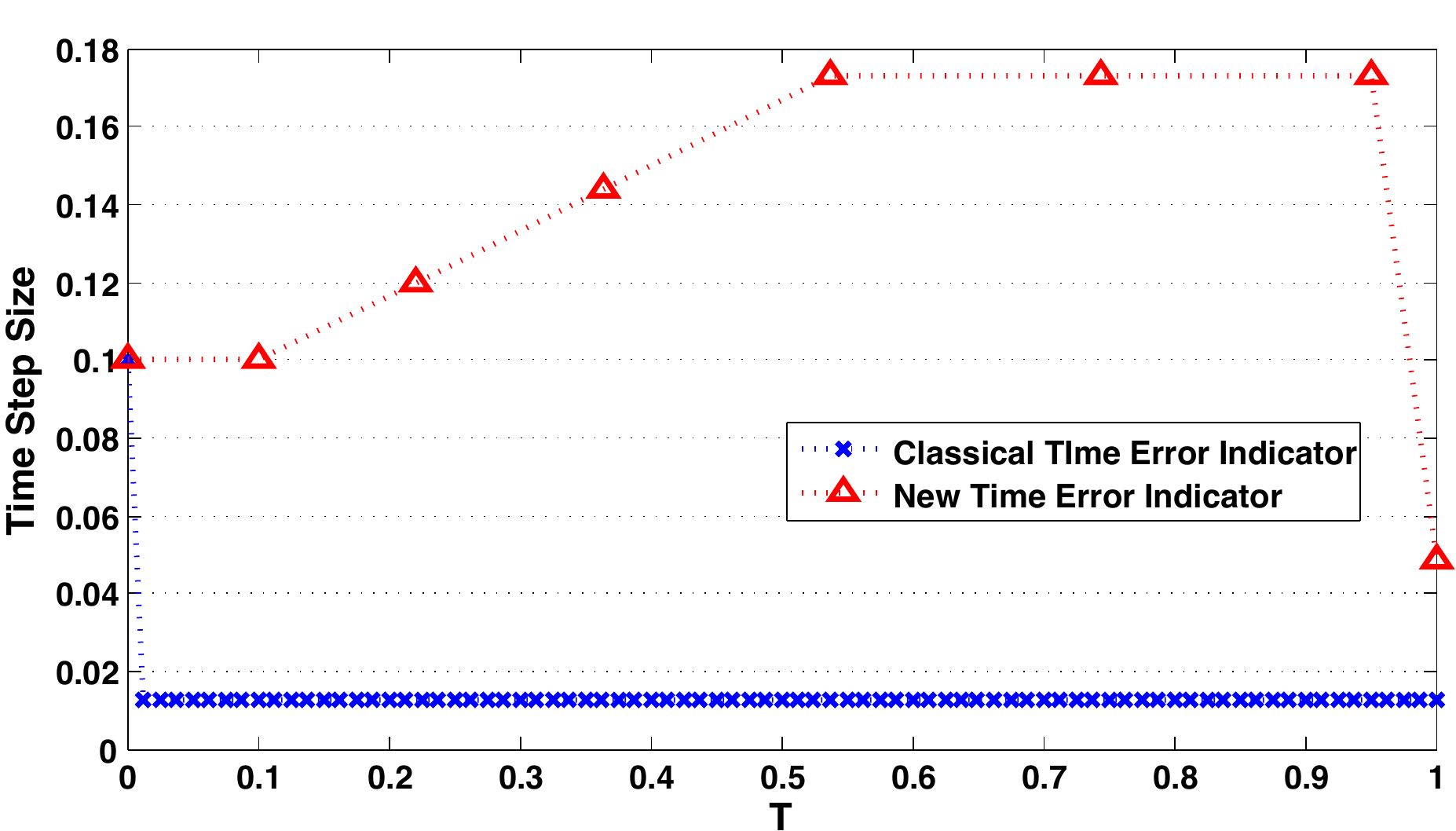}
\hfill
\includegraphics[width=0.49\textwidth, height=2.5in]{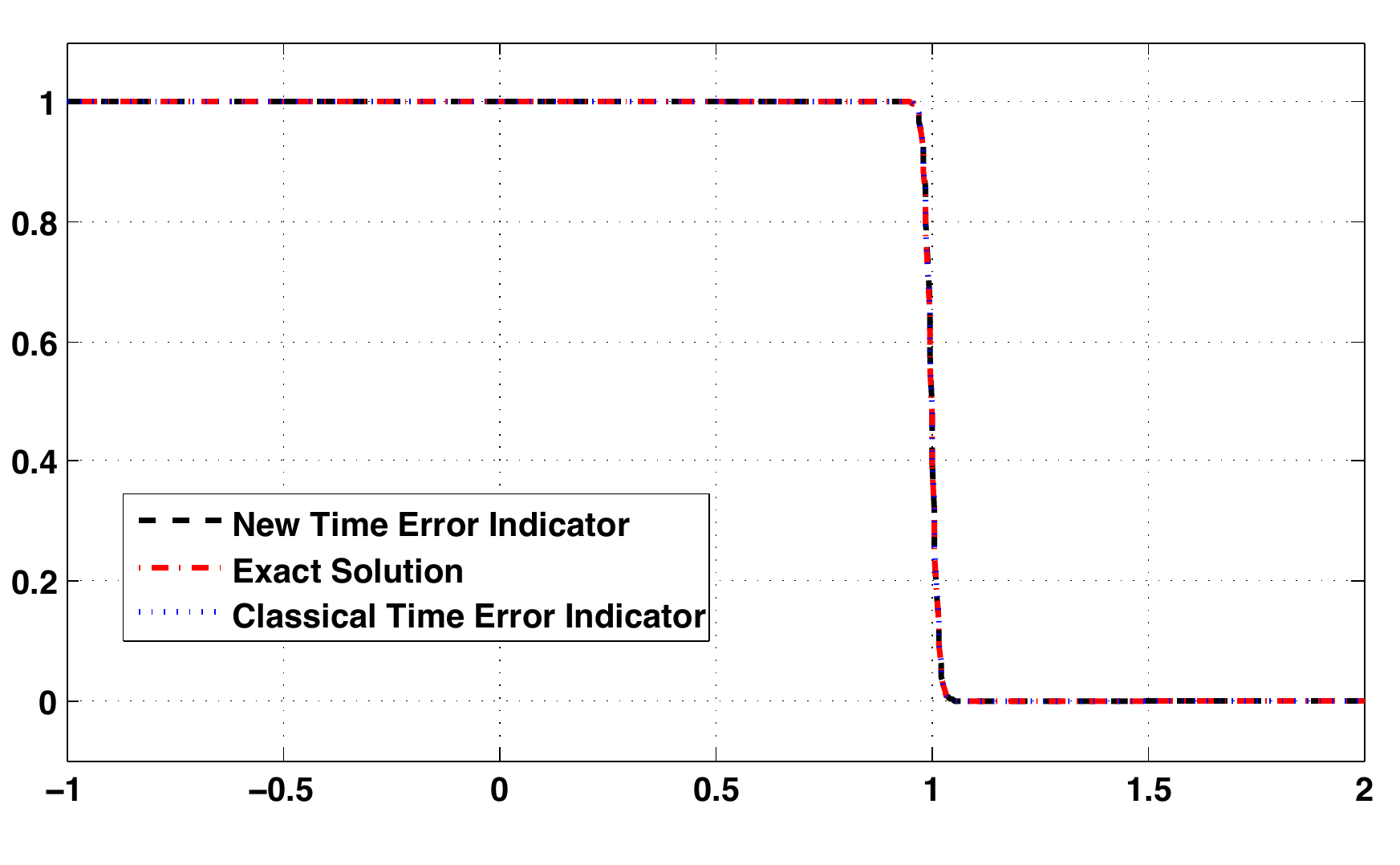}
\end{figure}

Figure \ref{fig:peak-1D} and \ref{fig:shock-1D} shows the numerical results for Example \ref{exp:peak-1D} and \ref{exp:shock-1D} respectively.  We can easily see that the new time error indicator allows larger time step size, but the numerical solution still maintains the accuracy.  This means the the new time error indicator is reliable and effective.

\subsubsection{2D Examples}
\begin{example}\label{exp:peak-2D}
This example is the benchmark Gaussian-cone problem (cf.~\cite{Xiu2001}).
Given the velocity field $b = (y, -x)^{T}$
and the initial condition 
\begin{equation*}
u(x,y,0) = \exp\left\{-[(x-x_{0})^{2}+(y-y_{0})^{2}]/(2\lambda^{2})\right\},
\end{equation*}
and $f=0$.  The exact solution is 
\begin{equation*}
u(x,y,t) = \frac{\lambda^{2}}{\lambda^{2}+2\varepsilon t}\exp\{-[\hat{x}^{2} + \hat{y}^{2}]/(2\lambda^{2}+4\varepsilon t)\},
\end{equation*}
where $\hat{x} = x-x_{0}\cos(t)-y_{0}\sin(t)$, $\hat{y} = y + x_{0}\sin(t)-y_{0}\cos(t)$, $\lambda = \frac{1}{8}$, and $(x_{0},y_{0}) = (-\frac{1}{2},0)$. The computational domain is $[0,1]\times[0,1]$ in $\mathbb{R}^{2}$.
\end{example}

The following two tests are generalizations of the one-dimensional problem, Problem~\ref{exp:shock-1D}. 
\begin{example}\label{exp:shock1-2D}
Given the velocity field $b = (1,0)^{T}$, and the initial condition
\begin{equation*}
u(x,y,0) = \left\{\begin{array}{rl}
1 \quad &\text{if } x <0.2\\
0 \quad &\text{otherwise},
\end{array}
\right.
\end{equation*}
and $f=0$. Then the exact solution is 
\begin{equation*}
\small
u(x,y,t) = \frac{1}{2}\left\{  \mathrm{erfc}\left(\frac{x-t}{2\sqrt{\varepsilon t}}\right) + \exp\left(\frac{x}{\varepsilon}\right)\mathrm{erfc}\left(\frac{x+t}{2\sqrt{\varepsilon t}}\right) \right\},
\end{equation*}
the computational domain is $[0,1]\times[0,1]$ in $\mathbb{R}^{2}$.
\end{example}

\begin{example}\label{exp:shock2-2D}
Given the velocity field $b = (1,1)^{T}$, and the initial condition as
\begin{equation*}
u(x,y,0) = \left\{\begin{array}{rl}
1 \quad &\text{if } x <0.2  \text{ and } y<0.2\\
0 \quad &\text{otherwise}.
\end{array}
\right.
\end{equation*}
and $f=0$. The exact solution is 
\begin{equation*}\footnotesize
u(x,y,t) = \frac{1}{4}\left\{  \mathrm{erfc}\left(\frac{x-t}{2\sqrt{\varepsilon t}}\right) + \exp\left(\frac{x}{\varepsilon}\right)\mathrm{erfc}\left(\frac{x+t}{2\sqrt{\varepsilon t}}\right) \right\}\left\{  \mathrm{erfc}\left(\frac{y-t}{2\sqrt{\varepsilon t}}\right) + \exp\left(\frac{y}{\varepsilon}\right)\mathrm{erfc}\left(\frac{y+t}{2\sqrt{\varepsilon t}}\right) \right\},
\end{equation*}
and the computational domain is $[0,1]\times[0,1]$ in $\mathbb{R}^{2}$.
\end{example}

\begin{figure}[ht]
\caption{Example \ref{exp:peak-2D} ($\varepsilon = 10^{-6}$). Left: time step size comparison; Middle: numerical solution; Right: adaptive mesh.} \label{fig:peak-2D}
\centering
\includegraphics[width=0.33\textwidth,height=2in]{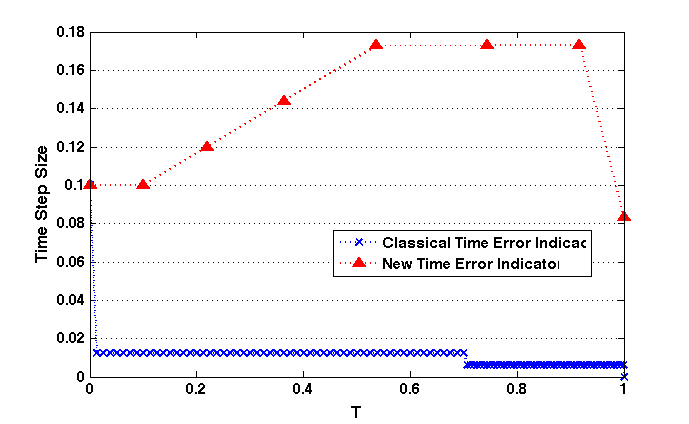}
\includegraphics[width=0.33\textwidth,height=2in]{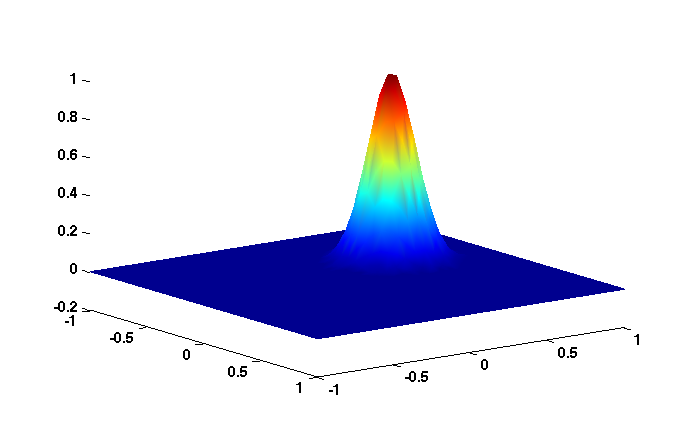}
\includegraphics[width=0.32\textwidth,height=2in]{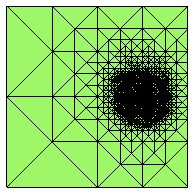}
\end{figure}

\begin{figure}[ht]
\caption{Example \ref{exp:shock1-2D} ($\varepsilon = 10^{-6}$). Left: time step size comparison; Middle: numerical solution; Right: adaptive mesh.} \label{fig:shock1-2D}
\centering
\includegraphics[width=0.33\textwidth, height=2in]{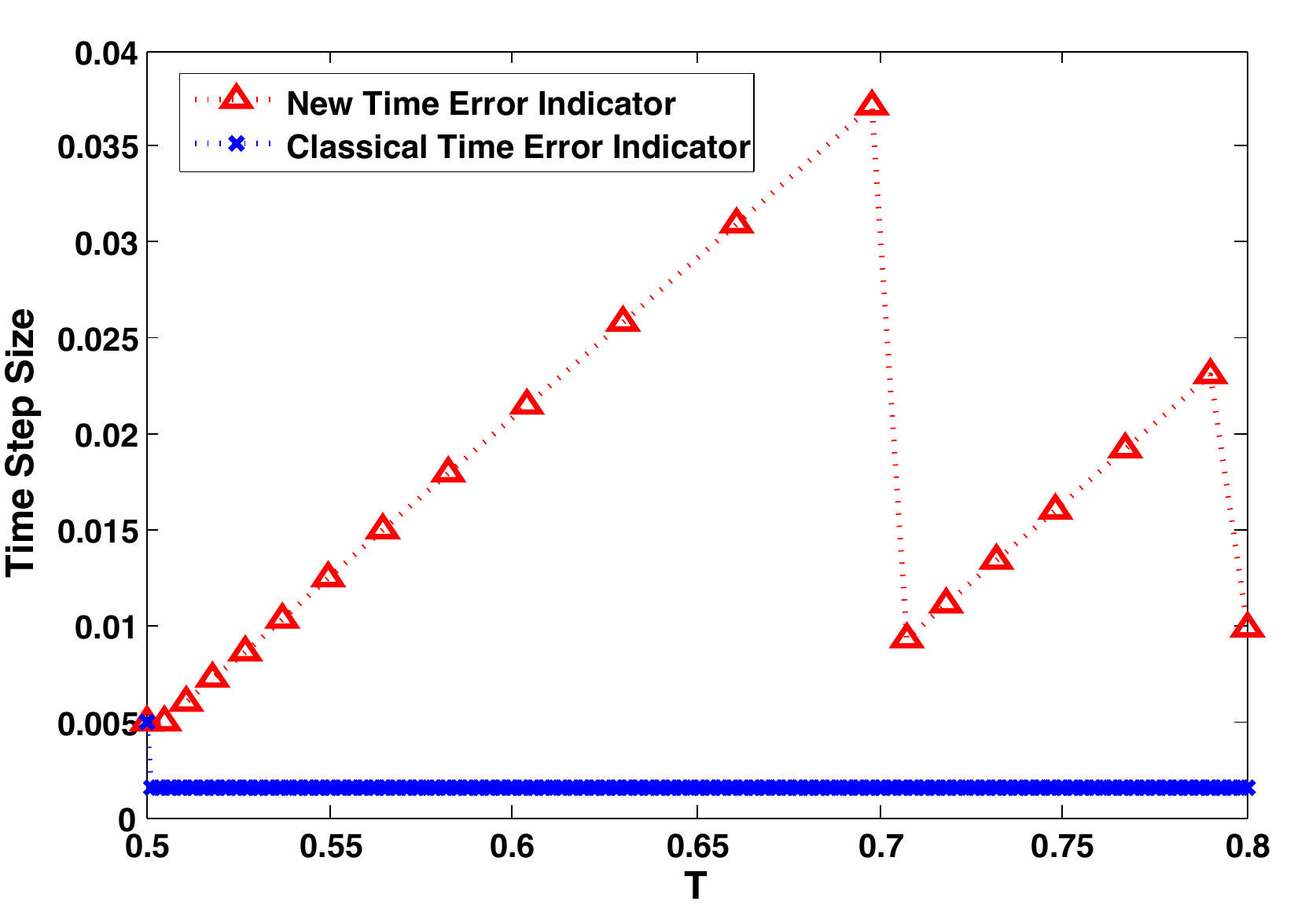}
\includegraphics[width=0.33\textwidth, height=2in]{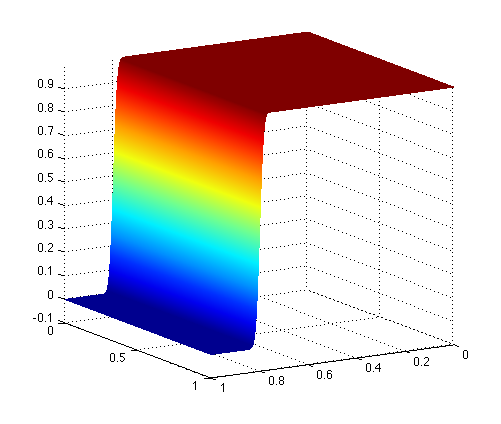}
\includegraphics[width=0.32\textwidth, height=2in]{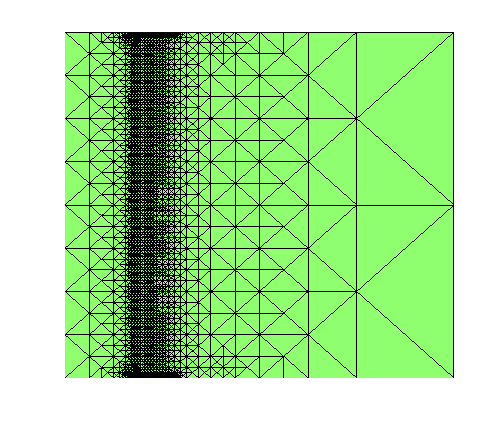}
\end{figure}

\begin{figure}[ht]
\caption{Example \ref{exp:shock2-2D} ($\varepsilon = 10^{-6}$). Left: time step size comparison; Middle: numerical solution; Right: adaptive mesh.}\label{fig:shock2-2D}
\centering
\includegraphics[width=0.33\textwidth, height=2in]{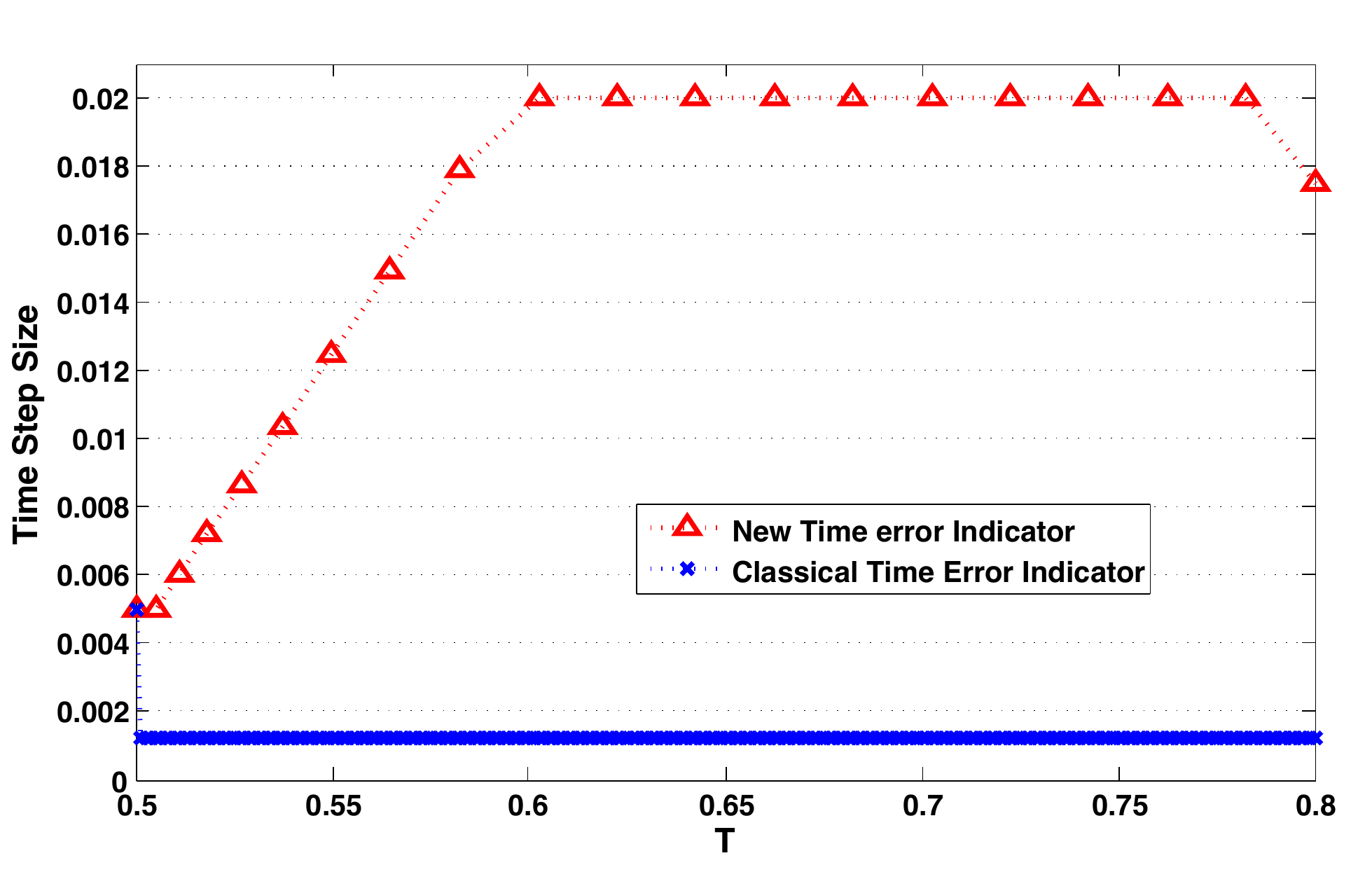}
\includegraphics[width=0.33\textwidth, height=2in]{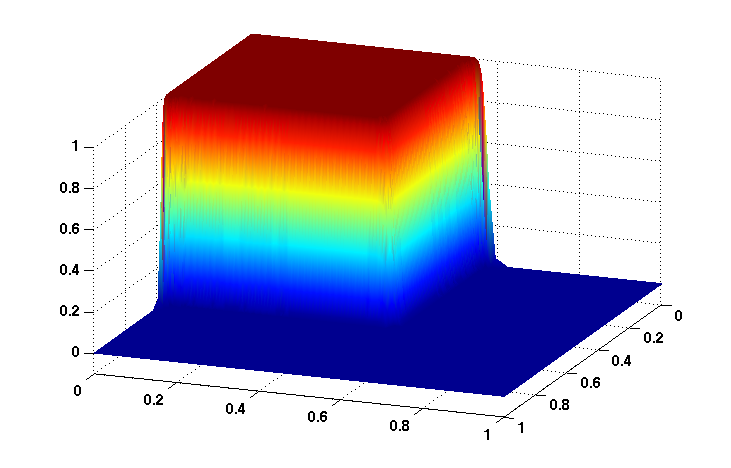}
\includegraphics[width=0.32\textwidth, height=2in]{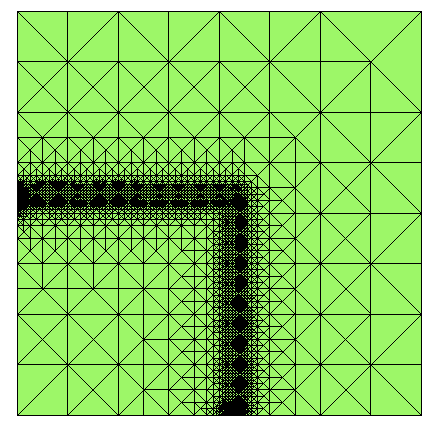}
\end{figure}

Figure \ref{fig:peak-2D}, \ref{fig:shock1-2D} and \ref{fig:shock2-2D} show the numerical results for Example \ref{exp:peak-2D}, \ref{exp:shock1-2D} and \ref{exp:shock2-2D} respectively.  As shown in the pictures, the new time error indicator allows larger time step size and the space error indicator captures the singularity.  Overall, the adaptive finite element method based on our error estimate is effective and reliable for convection dominated diffusion problems.

\section{Conclusion Remarks}

In this paper, we discuss the adaptive ELM for linear convection-diffusion equations.  We 

\begin{itemize}
\item derive a new temporal error indicator along the characteristics. We are able to show optimal convergence rate of temporal semi-discretization with minimal regularity of the exact solution.

\item combine the new temporal error indicator with residual-type spatial error indicator and obtain a posteriori error estimation for the fully discretization ELM.  

\item design efficient adaptive algorithm based on the a posteriori error estimators.  Numerical results shows robustness of the new adaptive method and allows larger time steps compared with standard temporal error indicator for ELM in general.

\end{itemize}

For the future work, we are working on generalize the a posteriori error estimation for nonlinear convection dominated problems, where the characteristics has to be solved approximately.  In this case, ODE solvers that preserves the determinant of the Jacobian of the flow will be important.  Meanwhile, we will also generalize the algorithm to Navier-Stokes equations. 

\smallskip
\noindent{\bf Acknowledgment.}
The authors would like to thank Professor Ricardo H. Nochetto and Professor Long Chen for their comments on earlier versions of this paper. NSF fund.

\bibliographystyle{abbrv}

\begin{thebibliography}{10}

\bibitem{Abbott1966}
M.~B. Abbott.
\newblock {\em {An introduction to method of characteristics}}.
\newblock American Elsevier, 1966.

\bibitem{Achdou2000}
Y.~Achdou and J.-L. Guermond.
\newblock Convergence analysis of a finite element
  projection/{L}agrange-{G}alerkin method for the incompressible
  {N}avier-{S}tokes equations.
\newblock {\em SIAM J. Numer. Anal.}, 37(3):799--826, 2000.

\bibitem{Adjerid1999}
S.~Adjerid, B.~Belguendouz, and J.~E. Flaherty.
\newblock A posteriori finite element error estimation for diffusion problems.
\newblock {\em SIAM J. Sci. Comput.}, 21(2):728--746, 1999.

\bibitem{Ainsworth2000}
M.~Ainsworth and J.~T. Oden.
\newblock {\em A posteriori error estimation in finite element analysis}.
\newblock Pure and Applied Mathematics (New York). Wiley-Interscience [John
  Wiley \& Sons], New York, 2000.

\bibitem{Bernardi2004}
C.~Bernardi and R.~Verf{\"u}rth.
\newblock A posteriori error analysis of the fully discretized time-dependent
  {S}tokes equations.
\newblock {\em M2AN Math. Model. Numer. Anal.}, 38(3):437--455, 2004.

\bibitem{Binev2004a}
P.~Binev, W.~Dahmen, and R.~Devore.
\newblock {Adaptive finite element methods with convergence rates}.
\newblock {\em Numerische Mathematik}, 97:219--268, 2004.

\bibitem{Brooks1982a}
N.~Brooks and T.~J. Hughes.
\newblock {Streamline Upwind/Petrov-Galerkink formulations for convection
  dominated flows with particular emphasis on the incompressilbe Navier-Stokes
  equations}.
\newblock {\em Computer Methods in Applied Mechanics and Engineering},
  32:199----259, 1982.

\bibitem{Cascon}
J.~M. Cascon, C.~Kreuzer, R.~H. Nochetto, and K.~G. Siebert.
\newblock Quasi-optimal convergence rate for an adaptive finite element method.
\newblock {\em SIAM J. Numer. Anal.}, 46(5):2524--2550, 2008.

\bibitem{Cawood2000}
M.~E. Cawood, V.~J. Ervin, W.~J. Layton, and J.~M. Maubach.
\newblock Adaptive defect correction methods for convection dominated,
  convection diffusion problems.
\newblock {\em Journal of Computational and Applied Mathematics}, 116(1):1--21,
  Apr. 2000.

\bibitem{Chen2004a}
Z.~Chen and J.~Feng.
\newblock An adaptive finite element algorithm with reliable and efficient
  error control for linear parabolic problems.
\newblock {\em Math. Comp.}, 73(247):1167--1193, 2004.

\bibitem{Chen2004}
Z.~Chen and G.~Ji.
\newblock Adaptive computation for convection dominated diffusion problems.
\newblock {\em Sci. China Ser. A}, 47(suppl.):22--31, 2004.

\bibitem{Chen2006}
Z.~Chen and G.~Ji.
\newblock Sharp {$L\sp 1$} a posteriori error analysis for nonlinear
  convection-diffusion problems.
\newblock {\em Math. Comp.}, 75(253):43--71 (electronic), 2006.

\bibitem{Chen2000}
Z.~Chen, R.~H. Nochetto, and A.~Schmidt.
\newblock A characteristic galerkin method with adaptive error control for the
  continuous casting problem.
\newblock {\em Computer Methods in Applied Mechanics and Engineering},
  189(1):249--276, 2000.

\bibitem{Cockburn1998a}
B.~Cockburn and C.-W. Shu.
\newblock {The local discontinuous Galerkin method for time-dependent
  convection-diffusion systems}.
\newblock {\em SIAM J. Numer. Anal.}, 35(6):2440----2463, 1998.

\bibitem{Davis1976}
G.~d.~V. Davis and G.~D. Mallinson.
\newblock {An evaluation of upwind and central difference approximations by a
  study of recirculating flow}.
\newblock {\em Computers \& Fluids}, 4:29--43, 1976.

\bibitem{Demkowicz1986}
L.~Demkowicz and J.~T. Oden.
\newblock An adaptive characteristic petrov-galerkin finite element method for
  convection-dominated linear and nonlinear parabolic problems in two space
  variables.
\newblock {\em Computer Methods in Applied Mechanics and Engineering},
  55(1-2):63--87, Apr. 1986.

\bibitem{Dorfler1996}
W.~D\"orfler.
\newblock A convergent adaptive algorithm for {Poisson}'s equation.
\newblock {\em SIAM Journal on Numerical Analysis}, 33:1106--1124, 1996.

\bibitem{Douglas1982}
J.~Douglas, Jr. and T.~F. Russell.
\newblock Numerical methods for convection-dominated diffusion problems based
  on combining the method of characteristics with finite element or finite
  difference procedures.
\newblock {\em SIAM J. Numer. Anal.}, 19(5):871--885, 1982.

\bibitem{Eriksson1991}
K.~Eriksson and C.~Johnson.
\newblock Adaptive finite element methods for parabolic problems i: A linear
  model problem.
\newblock {\em SIAM Journal on Numerical Analysis}, 28(1):43--77, 1991.

\bibitem{Eriksson1995}
K.~Eriksson and C.~Johnson.
\newblock Adaptive finite element methods for parabolic problems ii: Optimal
  error estimates in $l_\infty l_2 $ and $l_\infty l_\infty $.
\newblock {\em SIAM Journal on Numerical Analysis}, 32(3):706--740, 1995.

\bibitem{Eriksson1995a}
K.~Eriksson and C.~Johnson.
\newblock Adaptive finite element methods for parabolic problems iv: Nonlinear
  problems.
\newblock {\em SIAM Journal on Numerical Analysis}, 32(6):1729--1749, 1995.

\bibitem{Etienne2006}
J.~Etienne, E.~Hinch, and L.~J.
\newblock A {L}agrangian-{E}ulerian approach for the numerical simulation of
  freesurface flow of a viscoelastic material.
\newblock {\em Journal of Non-Newtonian Fluid Mechanics}, 136:157--166, 2006.

\bibitem{Ewing1983}
R.~E. Ewing.
\newblock {\em The Mathematics of Reservoir Simulation}.
\newblock SIAM, Philadelphia, 1983.

\bibitem{Feng1995}
K.~Feng and Z.-j. Shang.
\newblock {Volume-preserving algorithms for source-free dynamical systems}.
\newblock {\em Numerische Mathematik}, 463:451--463, 1995.

\bibitem{Hebeker1999}
F.~K. Hebeker and R.~Rannacher.
\newblock An adaptive finite element method for unsteady convection-dominated
  flows with stiff source terms.
\newblock {\em SIAM Journal on Scientific Computing}, 21(3):799--818, 1999.

\bibitem{Houston2001}
P.~Houston and E.~S{\"u}li.
\newblock Adaptive {L}agrange-{G}alerkin methods for unsteady
  convection-diffusion problems.
\newblock {\em Math. Comp.}, 70(233):77--106, 2001.

\bibitem{Jakob1959}
M.~Jakob.
\newblock {\em Heat Transfer}.
\newblock Wiley, New York, 1959.

\bibitem{Jia2011}
J.~Jia, X.~Hu, J.~Xu, and C.-S. Zhang.
\newblock {Effects of integrations and adaptivity for the Eulerian --
  Lagrangian method}.
\newblock {\em Journal of Computational Mathematics}, 29:367----395, 2011.

\bibitem{Johnson1990}
C.~Johnson, Y.-Y. Nie, and V.~Thom\'{e}e.
\newblock An a posteriori error estimate and adaptive timestep control for a
  backward euler discretization of a parabolic problem.
\newblock {\em SIAM Journal on Numerical Analysis}, 27(2):277--291, 1990.

\bibitem{Karakatsani2007}
F.~Karakatsani and C.~Makridakis.
\newblock A posteriori estimates for approximations of time-dependent stokes
  equations.
\newblock {\em IMA Journal of Numerical Analysis}, 27(4):741--764, 2007.

\bibitem{Kharrat2005}
N.~Kharrat and Z.~Mghazli.
\newblock Residual error estimators for the time-dependent stokes equations.
\newblock {\em Comptes Rendus Mathematique}, 340(5):405--408, 2005.

\bibitem{Kruger2003}
O.~Kr{\"u}ger, M.~Picasso, and J.-F. Scheid.
\newblock A posteriori error estimates and adaptive finite elements for a
  nonlinear parabolic problem related to solidification.
\newblock {\em Comput. Methods Appl. Mech. Engrg.}, 192(5-6):535--558, 2003.

\bibitem{Lakkis2006}
O.~Lakkis and C.~Makridakis.
\newblock Elliptic reconstruction and a posteriori error estimates for fully
  discrete linear parabolic problems.
\newblock {\em Math. Comp.}, 75(256):1627--1658, 2006.

\bibitem{Lee2006}
Y.~Lee and J.~Xu.
\newblock New formulations, positivity preserving discretizations and stability
  analysis for non-newtonian flow models.
\newblock {\em Comput. Methods Appl. Mech. Engrg.}, 195:1180--1206, 2006.

\bibitem{Moore1994}
P.~K. Moore.
\newblock A posteriori error estimation with finite element semi- and fully
  discrete methods for nonlinear parabolic equations in one space dimension.
\newblock {\em SIAM Journal on Numerical Analysis}, 31(1):149--169, 1994.

\bibitem{Morin2002}
P.~Morin, R.~H. Nochetto, and K.~G. Siebert.
\newblock Convergence of adaptive finite element methods.
\newblock {\em SIAM Review}, 44(4):631--658, 2002.

\bibitem{Nochetto2000}
R.~H. Nochetto, G.~Savar{\'e}, and C.~Verdi.
\newblock A posteriori error estimates for variable time-step discretizations
  of nonlinear evolution equations.
\newblock {\em Comm. Pure Appl. Math.}, 53(5):525--589, 2000.
\newblock Printed.

\bibitem{Nochetto2009}
R.~H. Nochetto, K.~G. Siebert, and A.~Veeser.
\newblock {Theory of adaptive finite element methods: An introduction}.
\newblock In {\em Multiscale, Nonlinear and Adaptive Approximation}, pages
  409----542. Springer Berlin Heidelberg, 2009.

\bibitem{Phillips1999}
T.~N. Phillips and A.~J. Williams.
\newblock Viscoelastic flow through a planar contraction using a
  semi-lagrangian finite volume method.
\newblock {\em Journal of Non-Newtonian Fluid Mechanics}, 87 (2-3):215--246,
  1999.

\bibitem{Phillips2000}
T.~N. Phillips and A.~J. Williams.
\newblock A semi-{L}agrangian finite volume method for {N}ewtonian contraction
  flows.
\newblock {\em SIAM J. Sci. Comput.}, 22(6):2152--2177, 2000.

\bibitem{Picasso1998}
M.~Picasso.
\newblock Adaptive finite elements for a linear parabolic problem.
\newblock {\em Computer Methods in Applied Mechanics and Engineering},
  167(3-4):223--237, 1998.

\bibitem{Pironneau1982}
O.~Pironneau.
\newblock On the transport-diffusion algorithm and its applications to the
  {N}avier-{S}tokes equations.
\newblock {\em Numer. Math.}, 38(3):309--332, 1982.

\bibitem{Roos2008}
H.~Roos, M.~Stynes, and L.~Tobiska.
\newblock {\em {Robust numerical methods for singularly perturbed differential
  equations: convection-diffusion-reaction and flow problems}}.
\newblock Springer-Verlag, 2008.

\bibitem{Schmidt2005}
A.~Schmidt and K.~G. Siebert.
\newblock {\em Design of adaptive finite element software}, volume~42 of {\em
  Lecture Notes in Computational Science and Engineering}.
\newblock Springer-Verlag, Berlin, 2005.
\newblock The finite element toolbox ALBERTA.

\bibitem{Sorbents1996}
I.~Sorbents, K.~U. Totsche, P.~Knabner, and I.~Kgel-knabner.
\newblock The modeling of reactive solute transport with sorption to mobile and
  immobile sorbents, 1. experimental evidence and model development.
\newblock {\em Water Resources Research}, 32:1611--1622, 1996.

\bibitem{Stevenson2005}
R.~Stevenson.
\newblock An optimal adaptive finite element method.
\newblock {\em SIAM Journal on Numerical Analysis}, 42(5):2188--2217, 2005.

\bibitem{Verfurth1996}
R.~Verf\"{u}rth.
\newblock {\em A review of a posteriori error estimation and adaptive mesh
  refinement techniques}.
\newblock Wiley and Teubner, 1996.

\bibitem{Verfurth1998}
R.~Verf{\"u}rth.
\newblock A posteriori error estimates for nonlinear problems: {$L\sp
  r(0,T;W\sp {1,\rho}(\Omega))$}-error estimates for finite element
  discretizations of parabolic equations.
\newblock {\em Numer. Methods Partial Differential Equations}, 14(4):487--518,
  1998.

\bibitem{VERFURTH2005}
R.~VERFURTH.
\newblock Robust a posteriori error estimates for nonstationary
  convection-diffusion equations.
\newblock {\em SIAM J. NUMER. ANAL.}, 43:1783--1802, 2005.

\bibitem{Wang2000}
H.~Wang.
\newblock An optimal-order error estimate for an {ELLAM} scheme for
  two-dimensional linear advection-diffusion equations.
\newblock {\em SIAM J. Numer. Anal.}, 37(4):1338--1368 (electronic), 2000.

\bibitem{Wang2007}
H.~Wang and K.~Wang.
\newblock Uniform estimates for eulerian--lagrangian methods for singularly
  perturbed time-dependent problems.
\newblock {\em SIAM Journal on Numerical Analysis}, 45(3):1305--1329, 2007.

\bibitem{Xiu2001}
D.~Xiu and G.~E. Karniadakis.
\newblock A semi-{L}agrangian high-order method for {N}avier-{S}tokes
  equations.
\newblock {\em J. Comput. Phys.}, 172(2):658--684, 2001.

\bibitem{Yeh1995}
G.-T. Yeh.
\newblock An adaptive local grid refinement based on the exact peak capture and
  oscillation free scheme to solve transport equations.
\newblock {\em Computer \& fluid}, 24:293, 1995.

\end{thebibliography}

\end{document}